\documentclass[aap,citesort,MSNbibl,seceqn,dvips]{arximspdf}

\doi{10.1214/10-AAP749}
\volume{21}
\issue{5}
\pubyear{2011}
\firstpage{2016}
\lastpage{2049}

\makeatletter

\def\bptnote#1{}

\newtheorem{theorem}{Theorem}[section]
\newtheorem{prop}[theorem]{Proposition}
\newtheorem{lem}[theorem]{Lemma}
\newproclaim{rem}[theorem]{Remark}
\newtheorem{cor}[theorem]{Corollary}
\newproclaim{defn}[theorem]{Definition}
\newproclaim{exm}[theorem]{Example}
\newproclaim{assumm}{Assumption}

\newcommand{\defi}{\stackrel{\triangle}{=}}
\newcommand{\R}{\mathbb{R}}

\makeatother

\begin{document}
\begin{frontmatter}

\title{Discounted continuous-time constrained Markov decision processes
in Polish spaces\thanksref{T1}}
\runtitle{Discounted continuous-time Markov decision processes}

\thankstext{T1}{Supported by NSFC, GDUPS(2010) and GRF(450508) from the
Research Grant Council of \mbox{HKSAR}.}

\begin{aug}
\author[A]{\fnms{Xianping} \snm{Guo}\ead[label=e1]{mcsgxp@mail.sysu.edu.cn}} and
\author[B]{\fnms{Xinyuan} \snm{Song}\corref{}\ead[label=e2]{xysong@sta.cuhk.edu.hk}}
\runauthor{X. Guo and X. Song}
\affiliation{Zhongshan University and The Chinese University of Hong Kong}
\address[A]{School of Mathematics\\
\quad and Computational Science\\
Zhongshan University\\
Guangzhou 510275\\
P. R. China\\
\printead{e1}} 
\address[B]{Department of Statistics\\
Chinese University of Hong Kong\\
LSB114, Shatin, Hong Kong\\
P. R. China\\
\printead{e2}}
\end{aug}

\received{\smonth{6} \syear{2010}}
\revised{\smonth{10} \syear{2010}}

%
\begin{abstract}
This paper is devoted to studying constrained continuous-time Markov
decision processes (MDPs) in the class of randomized policies depending
on \textit{state histories}. The transition rates may be \textit{unbounded},
the reward and costs are admitted to be \textit{unbounded from above and
from below}, and the state and action spaces are Polish spaces. The
optimality criterion to be maximized is the expected discounted
rewards, and the constraints can be imposed on the expected discounted
costs. First, we give conditions for the nonexplosion of underlying
processes and the finiteness of the expected discounted rewards/costs.
Second, using a technique of occupation measures, we prove that the
constrained optimality of continuous-time MDPs can be transformed to an
\textit{equivalent} (optimality) problem over a class of probability
measures. Based on the equivalent problem and a so-called \textit{$\bar
w$-weak convergence} of probability measures developed in this paper,
we show the existence of a constrained optimal policy. Third, by
providing a linear programming formulation of the equivalent
problem, we show the solvability of constrained optimal policies.
Finally, we use two \textit{computable} examples to illustrate our main
results.
\end{abstract}

%
\begin{keyword}[class=AMS]
\kwd{90C40}
\kwd{60J27}.
\end{keyword}
\begin{keyword}
\kwd{Continuous-time Markov decision process}
\kwd{unbounded transition rates}
\kwd{occupation measure}
\kwd{linear programming formulation}
\kwd{constrained optimal policy}.
\end{keyword}

\end{frontmatter}

\section{Introduction}

\textit{Constrained} Markov decision processes (MDPs) form an
important class of stochastic control problems and have been
widely studied. Existing works on constrained MDPs can be roughly
classified into four groups: (i)~constrained discrete-time MDPs
with denumerable states \cite
{A99,AS91,CF07,F00,FS99,FS96,FS95,HS89,K06,P94,Sen,ZS06}
and their extensive references, (ii) constrained
discrete-time MDPs with a Polish state space \cite{OJ00,OJ03,KNH,Pi97}
and their bibliographies, (iii) constrained continuous-time
MDPs with denumerable states \cite{GH09,GH03c,Pi98,PH08,ZG08}, and
(iv) constrained continuous-time MDPs with a Polish state space
\cite{G07a}. A review of these references shows that most of the related
literature is concentrated with the first three groups. To the best
of our knowledge, the fourth group is addressed only in~\cite{G07a}
for the average criteria. Concerning group (i), the existence and
algorithms of constrained optimal policies are given in
\cite{CF07,F00,FS99,FS96,FS95} for \textit{variant} discounted criteria
when states and actions are finite, in \cite{A99,K06,P94} for the
discounted criteria and denumerable states, and in
\cite{A99,AS91,HS89,P94,Sen} for the average criteria and
denumerable states. Also, the existence of constrained optimal
policies and linear programming formulation for group (ii) are
given in \cite{OJ00,Pi97} for the discounted criteria and in
\cite{OJ03,KNH,Pi97} for the average criteria. Although group
(iii) has been studied in \cite{GH09,GH03c,Pi98,PH08,ZG08}, the
references \cite{GH09,GH03c,Pi98,PH08,ZG08} deal with the case of a
single constraint, the transition rates in \cite{Pi98} are assumed
to be bounded, and the assumption of denumerable states in these
references cannot be dropped. On the other hand, as mentioned above,
constrained MDPs in Polish spaces are also studied in \cite{OJ00,OJ03,KNH,Pi97}
for the discrete-time case and in \cite{G07a} for
the continuous-time case. However, the reward and cost functions in
\cite{KNH} are assumed to be all \textit{bounded}, and all cost
functions in \cite{G07a,OJ00,OJ03,Pi97} are assumed to be
essentially \textit{nonnegative}. Further, such nonnegativeness
assumption cannot be removed because it is required for the use of
the standard weak convergence of probability measures. This in turn
implies that the constrained optimality problem of minimizing
nonnegative costs in \mbox{\cite{G07a,OJ00,OJ03}} with constraints imposed
on other nonnegative costs cannot be transformed to an equivalent
optimality problem of maximizing \textit{bounded} rewards as in
\cite{KNH} with constraints imposed on bounded costs. Hence, the
constrained discrete and continuous time MDPs with Polish spaces, in
which rewards (to be maximized) and costs (with constraints) may be
unbounded from above and from below, have not been studied.

On the other hand, as is known, continuous-time MDPs in Polish
spaces have been studied in \cite{G07b,G07a,GR06,K85,Pi98}.
However, the treatments in \cite{G07b,GR06,K85} are on the
unconstrained case, whereas the results in \cite{G07a} for the
constrained case cannot be applied to the case in which the
criterion to be maximized is \textit{unbounded} rewards. This is
because the
cost to be minimized in \cite{G07a} is required to be \textit
{nonnegative}. Moreover, the study in \cite{G07b,G07a,GR06}
with unbounded transition rates is limited to the class of \textit
{Markov} policies, and yet the case of
randomized policies depending on state histories in \cite{K85,Pi98} is
for \textit{bounded} transition rates. Hence, as
noted in \cite{GH09,GHP06,Y77}, the study on unconstrained
continuous-time MDPs with
unbounded transition rates and history-dependent policies is an
unsolved problem.

Constrained continuous-time MDPs with unbounded
transition rates and policies depending on state histories have \textit{not}
been studied yet, and they will be considered in this paper. More
precisely, we will deal with constrained continuous-time MDPs,
which have the following features: $(1)$ the transition rates may
be \textit{unbounded}; $(2)$ the reward and costs are admitted to be
unbounded \textit{from above and from below}; $(3)$ the state and
action spaces are Polish spaces; $(4)$ admissible policies can
be randomized and \textit{depend on state histories}; and $(5)$ the optimality
criterion is to \textit{maximize} expected discounted rewards, and
several constraints are imposed on expected discounted costs.

First, we give the conditions under which we ensure the nonexplosion
of underlying processes induced from \textit{unbounded} transition rates
and randomized policies \textit{depending on
state histories} (see Theorem \ref{t2} below). This result is a natural
extension
of the corresponding regularity of a jump Markov process in
\cite{C04,G07b,GH09,GR06,LMT96} to a so-called ``non-Markov'' case
and also a generalization of the regularity in
\cite{havis,K71,K85,KR95,LP01,Pi98,P94,S99,Y77} for \textit{bounded}
transition rates. Inspired by the condition for the nonexplosion, we
obtain a condition (see Theorem \ref{t3} below) for the
finiteness of the expected discount rewards/costs of each policy when
rewards/costs are \textit{unbounded}.

Second, as in \cite{A99,AS91,OJ00,OJ03,OL96,KNH,Pi05,Pi97} for
constrained MDPs, by introducing an occupation measure, we
prove that the constrained optimality problem in continuous-time MDPs
[see (\ref{P}) below] can be transformed into an \textit{equivalent}
optimality problem [see (\ref{4.2b}) below] over a class of some
probability measures.
The standard weak convergence technique used
in \cite{G07a,OJ00,OJ03,OL99,K85,KNH} for nonnegative costs does not
apply directly to the case wherein rewards/costs are unbounded from
above and from below. Therefore, to solve the equivalent optimality
problem in which rewards/costs may be unbounded from above and from
below, we introduce (Definition \ref{De-4.2} below) a
so-called \textit{$\bar w$-weak convergence} of probability measures.
This $\bar w$-weak convergence is an
extension of the standard weak convergence of probability measures.
Using the properties of the $\bar w$-weak convergence and occupation
measures developed here (see Theorem \ref{Th-4.1} and Lemmas \ref
{Th-4.5} and \ref{Th-4.4} below), we prove the existence of a
constrained optimal policy under mild reasonable conditions (see
Theorem \ref{Th-4.6} below). These conditions are slightly different
from the usual continuity-compactness ones in
\cite{G07b,GH09,GH03a,GH03c} for continuous-time MDPs and in \cite
{A99,AS91,OJ00,OJ03,OL99,KNH} for the discrete-time MDPs, and thus they
are weaker than those
in the literature \cite{G07b,GH09,GH03a,GH03c,P94}; see Remarks \ref
{rem5.1-new} and \ref{rem5.1} for details.

Third, for the solvability of constrained optimal policies, we further
transform the equivalent optimality problem to a linear programming
(LP) problem [see (\ref{LP}) below] by using the properties of
occupation measures
again. Then we present the relationship between a constrained optimal
policy and an optimal solution to the LP (see Theorem \ref{Th-5.1}
below), and characterize a stationary policy (see Theorem~\ref{Th-3.2}
below). This
relationship and characterization of a stationary policy are used to
obtain the solvability and structure of a constrained optimal policy
(see Corollary \ref{Co-1} and Theorem \ref{Th-5.2}
below).

Finally, to illustrate our main results, we present two \textit
{computable} examples in which our conditions are
satisfied, whereas some of those in \cite{G07a,OJ00,OJ03,OL99,K85,KNH}
fail to hold (see Remark \ref{rem5.2} below).
In particular, our approach is also \textit{suitable} to the case of
discrete-time MDPs with rewards/costs
being unbounded from above and from below, and similar results for the
discrete-time case can also be obtained; see Remark \ref{rem5.1-rr}
for details.
However, our model \textit{cannot} be transformed to an equivalent one
of discrete-time MDPs
using the uniformization technique because the transition rates in our
model may be \textit{unbounded}.

The rest of this paper is organized as follows. In Section \ref
{section-2}, the
model and the constrained optimality
problem that we are concerned with are introduced. The main results of
this paper are stated in Section \ref{section-3}, and illustrated with
computable
examples in Section \ref{section-4}. The proofs of the main results are
presented in
Section \ref{proof}.

\section{The model for constrained continuous-time MDPs}\label{section-2}

\mbox{}

\textit{Notation.} If ${\mathrm X}$ is a Polish space (i.e., a complete
and separable
metric space) and $\bar w\geq1$ is a real-valued measurable function
on $X$, we denote
by ${\mathcal{B}}({\mathrm X})$ the Borel $\sigma$-algebra on $X$,
by $D^c$ the complement of a set $D\subseteq X$ (with respect to~$X$),
by $\|u\|_{\bar w}$
the $\bar w$-weighted norm of a real-valued measurable function
$u$ on $X$ [i.e., $\|u\|_{\bar w}:={\sup_{x\in X}}|u(x)|/\bar
w(x)$], by $C_b(X)$ the set of all bounded
continuous functions on $X$, and by ${\mathcal{P}}(X)$ the set of all
probability measures
on ${\mathcal{B}}({\mathrm X})$. Let
\[
B_{\bar
w}(X):=\{u| \|u\|_{\bar w}<\infty\}
\]
be the Banach space.

We now introduce the model of constrained continuous-time MDPs,
%
%
\begin{equation} \label{2.1}\quad
\bigl\{S, \bigl(A(x)\subseteq A, x\in S\bigr), q(\cdot|x,a), r(x,a),
\bigl(c_n(x,a),d_n, 1\leq n\leq N\bigr) \bigr\},
\end{equation}
where $S$ is a \textit{state space}, $A$ is an \textit{action space}, and
$A(x)$ is a Borel \textit{set of admissible actions} at state $x\in
S$. We suppose that $S$ and $A$ are Polish spaces, and the
following
set:
%
%
\begin{equation}\label{2.2}
K:=\{(x,a)| x\in S, a\in A(x)\}
\end{equation}
is a Borel subset of $S\times A$.

The function $q(\cdot|x,a)$ in (\ref{2.1}) refers to \textit{transition
rates}, that is, it satisfies the following:
\begin{longlist}[${(T_3)}$]
\item[${(T_1)}$] For each fixed $(x,a)\in K, q(\cdot|x,a)$ is a
\textit{signed}
measure on ${\mathcal{B}}(S)$, whereas for each fixed
$D\in{\mathcal{B}}(S)$, $ q(D|\cdot)$ is a real-valued
Borel-measurable function on $K$;
\item[${(T_2)}$] $0\leq
q(D|x,a)<\infty$ for all $(x,a)\in K$ and $ x\notin D \in
{\mathcal{B}}(S)$; and
\item[${(T_3)}$] $q(S|x,a)=0 $ for all
$(x,a)\in K$. [Hence, $q(\{x\}|x,a)$ is finite for all $(x,a)\in
K$.]
\end{longlist}
The model is also assumed to be \textit{stable}, which means
%
%
\begin{equation}\label{2.3}
q^*(x):={\sup_{a \in A(x)}}|q(\{x\}|x,a)|<\infty\qquad\forall
x\in S.
\end{equation}
Finally, the function $r(x,a)$ on $K$ denotes the reward, whereas
the functions $c_n(x,a)$ on $K$ and the real
numbers $d_n$ denote the costs and constraints, respectively. We
assume that $r(x,a)$ and $c_n(x,a)$ are real-valued measurable on $K$.
[\textit{$r(x,a)$ is allowed to take positive and negative values},
\textit{so it
can be interpreted as a} \textit{cost} \textit{rather than a} ``\textit
{reward}'' \textit{only}.]

To complete the specification of the constrained optimality problem,
we of course need an optimality criterion. This requires the
definition of a class of policies admissible to a controller. To do
so, we introduce some notation as in \cite{b0,K85,KR95}.

Let $S_\infty:=S\cup\{x_\infty\}$ with $x_\infty$ being an isolated
point, $\Omega^0:=(S\times\R_+)^\infty$ with
$\R_+:=(0,\infty)$ and $\Omega:=\Omega^0 \cup\{(x_0,\theta
_1,x_1,\ldots,\theta_{k-1},x_{k-1},\infty,x_\infty,
\ldots)| \theta_l\in\R_+, x_0, x_l\in S$ for each $1\leq l\leq k-1$
and $k\geq2\}$. By the corresponding modification of the
$\sigma$-algebra over $\Omega^0$, we can obtain the basic
measurable space $(\Omega,\mathcal{F})$. Then we define maps $T_k,
X_k, \Theta_k$ $(k=0,1,\ldots)$ and $\xi_t$ $(t\geq0)$ on
$(\Omega,\mathcal{F})$ as follows: for each $e:=(x_0,\theta_1, x_1,
\ldots, \theta_k,x_k, \ldots)\in\Omega$, let
%
%
\begin{eqnarray}
T_k(e)&:=&\theta_1+\cdots+\theta_k \qquad\mbox{(for
$k\geq1$)},\nonumber\\[-8pt]\\[-8pt]
T_\infty(e)&:=&\lim_{k\to\infty} T_k(e)\qquad \mbox{with }
T_0(e):=0; \nonumber\\
X_{k-1}(e)&:=&x_{k-1},\qquad \Theta_k(e):=\theta_k\qquad
\mbox{for } k\geq1; \nonumber\\
\label{2.4}
\xi_t(e)&:=&\sum_{k\ge0} x_kI_{\{T_k\le t<T_{k+1}\}}(e)+x_\infty
I_{\{T_\infty\le t\}}(e),
\end{eqnarray}
where $I_D$ stands for the indicator function of a set $D$. Let
$h_k(e)=(x_0,\theta_1,x_1,\ldots,\break\theta_{k},x_{k})$, and call
$h_k(e)$ a
$k$-component \textit{state history}. Obviously, these maps are
measurable on
${\mathcal{F}}$. In what follows, the argument $e=(x_0,\theta_1, x_1,
\ldots,\theta_k,x_k,\ldots)$ is often omitted.

Components $\Theta_k$ play the role of inter-jump intervals or
sojourn times, $T_k$ are the jump epoches, and $X_k$ denotes the
state of the process $\{\xi_t,t\geq0\}$ on $[T_k,T_{k+1})$. We do
not intend to consider the process after moment $T_\infty$, so we view
it to be absorbed in state $x_\infty$. Hence, we write
$q(\cdot|x_\infty,a_\infty)\equiv0$, where $a_\infty$ is an
isolated point, and let $A(x_\infty):=\{a_\infty\}$,
$A_\infty:=A\cup\{a_\infty\}$.

Let $\R_+^0:=[0,\infty)$, and introduce the integer-valued random
measure $\mu^*$ on $\R_+^0\times S$ by
%
%
\begin{equation}\label{2.5}
\mu^*(dt,dx)=\sum_{k\ge0}I_{\{T_k<\infty\}}\delta_{(T_k,X_k)}(dt,dx),
\end{equation}
where $\delta_y(\cdot)$ is the Dirac measure concentrated at any
point $y$. Then we take the right-continuous family of
$\sigma$-algebras $\{\mathcal{F}_t\}_{t\ge0}$ with $\mathcal
{F}_t:=\sigma\{\mu^*([0,s]\times D),s\in[0,t],D\in{\mathcal
{B}}(S)\}$,
and let
\[
{\mathcal{P}}:=\sigma\bigl(B\times\{0\}, C\times(s,\infty)| B\in
{\mathcal{F}}_0, C\in{\mathcal{F}}_{s-}, s>0\bigr),
\]
where
${\mathcal{F}}_{s-}:=\bigvee_{t<s}{\mathcal{F}}_t$. Then, as in
\cite{b0,K85,KR95}, a real-valued function on $\Omega\times
\R^0_+$ is called \textit{predictable} if it is measurable with
respect to $\mathcal P$.

We next introduce the definition of a policy, which is the same as
in \cite{K85} and a generalization of the corresponding one in
\cite{KR95,Pi05,Pi98} for denumerable states.
\begin{defn}\label{d1}
A transition probability $\pi$ from
$(\Omega\times\R^0_+,\mathcal{P} )$ onto $(A_\infty,\break
{\mathcal{B}}(A_\infty))$ such that $\pi(A(\xi_{t-}(e))|e,t)\equiv
1$ is
called a policy, which can be randomized and depend on state histories.
A policy is
called \textit{randomized stationary} if there exists a transition
probability $\phi$ from $(S,{\mathcal{B}}(S))$ onto $(A, {\mathcal{B}}(A))$
such that $\phi(A(x)|x)\equiv1$ and
$\pi(da|e,t)=I_{\{t<T_\infty\}}(e)\phi(da|\xi_{t-}(e))+I_{\{t\geq
T_\infty\}}(e)\delta_{a_\infty}(da)$. We will write such a
randomized stationary policy as $\phi$.
A randomized stationary policy $\phi$ is called (deterministic) \textit
{stationary} if there exists a measurable
function $f$ from $(S,{\mathcal{B}}(S))$ onto $(A, {\mathcal{B}}(A))$
such that $\phi(\{f(x)\}|x)\equiv1$. Such
a stationary policy will be written as $f$.
\end{defn}

We denote by $\Pi, \Pi_s$ and $F$ the classes of all policies,
randomized stationary policies and stationary policies,
respectively. Equivalently, $\Pi_s$ is the set of all stochastic
kernels $\phi$ on $A$ given $S$ such that $\phi(A(x)|x)=1$ for all
$x\in S$, and $F$ is the set of all measurable functions $f$ from
$S$ to $A$ such that $f(x)\in A(x)$ for all $x\in S$. Obviously,
$F \subset\Pi_s \subset\Pi$.
\begin{rem}\label{r2}
The requirement of predictability of a
policy implies that at time $t\geq0$ each policy depends on only
the past jump moments $T_0,T_1,\ldots,\break T_m \leq t$ and the
corresponding states $x_0,\ldots,x_m\in S$. This means that a policy
may depend on state histories. However, the class $\Pi$ is not the complete
collection of all history-dependent policies. This is because
each state history
$h_k=(x_0,\theta_1,x_1,\ldots,\theta_{k},x_{k})$ does not include
past actions $a_m$ $(0\leq m\leq k)$. To overcome the shortcoming
of the definition of a state history, a possible and natural way is to
replace $h_k$ with a new history
$(x_0,a_0,\theta_1,\ldots,x_{k-1},a_{k-1},\theta_{k},x_{k})$
including past actions. If we do so, some results in \cite{b0,KR95}
such as the structure of the probability measure $P^\pi_\gamma$ in
(\ref
{3.2}) and the predictable properties of the randomized measure $\nu
^\pi
$ in (\ref{2.6}) and functions $m(D|e,t)$ in (\ref{3.1}), which are
required in following arguments, need to be checked one by one. Since
these desired results
for the case of new histories have
not been proven, we still use the definition of a policy in Definition
\ref{d1}, which is the same as in \cite{K85,KR95,Pi05,Pi98}, and
which is also a generalization of the corresponding one in
\cite{C04,G07b,G07a,GH09,GHP06} for a Markov policy.
\end{rem}

For each $\pi\in\Pi$, by Definition
\ref{d1} we see that the random measure on
$\R^0_+\times S$ given by
%
%
\begin{eqnarray}\label{2.6}
\nu^\pi(e,dt,D):=\biggl[\int_A\pi(da|e,t)q(D|\xi_{t-}(e),a)I_{\{
\xi
_{t-}\notin
D\}}(e) \biggr]\,dt\nonumber\\[-8pt]\\[-8pt]
&&\eqntext{\mbox{for } D\in{\mathcal{B}}(S)}
\end{eqnarray}
is predictable, and $\nu^\pi(\{t\}\times
S)=\nu^\pi([T_\infty,\infty)\times S)\equiv0$ for all $t\geq0$.
Thus, for any initial distribution $\gamma\in{\mathcal{P}}(S)$,
Theorem 4.27 in \cite{KR95} (or Theorem 3.6 in \cite{b0}) ensures
the existence of a unique probability measure $P_{\gamma}^\pi$ on
$(\Omega,\mathcal{F})$ such that $P_{\gamma}^\pi\{x_0\in
dx\}=\gamma(dx)$, and $\nu^\pi$ is a dual predictable projection
of the measure $\mu^*$ in (\ref{2.5}). The expectation operator
with respect to $P_\gamma^\pi$ is denoted by
$E_\gamma^\pi$. In particular, $E_\gamma^\pi$ and
$P_\gamma^\pi$ will be written as $E_x^\pi$ and $P_x^\pi$,
respectively, when $\gamma$ is the Dirac measure located at point
$x\in S$.

For any fixed\vspace*{1pt} $\pi\in\Pi$ and $\gamma\in{\mathcal{P}}(S)$, let us
recall how the measure $P^\pi_\gamma$ is constructed. First, by
Definition \ref{d1} we see that, for each fixed $D\in
{\mathcal{B}}(S)$, the following function on $\Omega\times\R^0_+$:
\[
m(D|e,t):=\int_A\pi(da|e,t)q(D|\xi_{t-}(e),a)I_{\{\xi_{t-}\notin
D\}}(e)
\]
is predictable, and thus (by Lemma 3.3 in \cite{b0}) has the
following representation:
%
%
\begin{eqnarray}\label{3.1}
m(D|e,t)&=:&I_{\{0\}}(t)m_0(D|x_0,0)\nonumber\\[-8pt]\\[-8pt]
&&{}+\sum_{k=0}^\infty I_{\{T_k<t\le
T_{k+1}\}}(e)m_k\bigl(D|h_k(e),t-T_k\bigr),\nonumber
\end{eqnarray}
where $m_k(\cdot|h_k(e),\tilde t\,)$ (depending on $\pi$) is a
measure on ${\mathcal{B}}(S)$ [for any fixed $h_k(e)$ and $\tilde t$],
$m_k(D|h_k(e),\tilde t\,)$ is measurable in $(e,\tilde t\,)$ [for
any fixed $D\in{\mathcal{B}}(S)$] and $m_k(\{x_k\}|h_k(e),\tilde
t\,)=0$ for all $x_k\in S$ and $k\geq0$. Let\vspace*{1pt} $\hat H_0\defi S, \hat
H_k\defi S\times(\R_+\times S_\infty)^k$ for $k\geq1$. Noting
that a\vspace*{2pt} measure $\gamma$ on ${\mathcal{B}}(\hat H_0)$ is given, we
suppose that the measure $P^\pi_{\gamma}$ on ${\mathcal{B}}(\hat H_k)$
has been constructed, then $P^\pi_{\gamma}$ on ${\mathcal{B}}(\hat
H_{k+1})$ is determined as follows:
%
%
\begin{eqnarray}\label{3.2}
&& P^\pi_{\gamma}\bigl(\Gamma\times(d\tilde t,dx)\bigr)\nonumber\\
&&\qquad:=\int_{\Gamma}
P^\pi_{\gamma}(dh_k)I_{\{\theta_{k+1}<\infty\}}m_k(dx|h_k,\tilde
t\,)e^{-\int_0^{\tilde t}m_k(S|h_k,v)\,dv}\,d\tilde
t;\nonumber\\[-8pt]\\[-8pt]
&&P^\pi_{\gamma}\bigl(\Gamma\times(\infty,x_\infty)\bigr)\nonumber\\
&&\qquad:=\int_{\Gamma
}P^\pi
_{\gamma}(dh_k)\{I_{\{\theta_{k+1}=\infty\}}+I_{\{\theta
_{k+1}<\infty\}}
e^{-\int_0^\infty m_k(S|h_k,v)\,dv}\}, \nonumber
\end{eqnarray}
where $\Gamma\in\mathcal{B}({\hat H}_k)$. According to the
Ionescu Tulcea theorem in \cite{b7}, there exists a unique
probability measure $P^\pi_{\gamma}$ on $(\Omega,\mathcal{F})$,
which has projections onto the spaces of $k$-component
state histories satisfying relations (\ref{3.2}).

For any given $\gamma\in{\mathcal{P}}(S)$ and $\pi\in\Pi$, using
(\ref{3.1}) and (\ref{3.2}), we now give a somewhat informal
description of how the process $\{\xi_t,t\geq0\}$ evolves. Suppose
that the process is at state $x_k$ at time $t\in[T_k,T_{k+1})$ $
(k\geq0)$. Then, a transition from $x_k$ to a set $D$ of states
occurs with probability $m_k(D|h_k,t-T_k)$, or the process remains
at $x_k$ with probability $1-m_k(S|h_k,t-T_k)\,dt+o(dt)$. In the
former case, the sojourn time $\Theta_{k+1}$ of $\{\xi_t,t\geq0\}$
at $x_k$ has a distribution with a so-called ``density function''
$e^{-\int_0^{t}m_k(S|h_k,v)\,dv}$.

As mentioned above, we do not intend to consider the process after
moment~$T_\infty$. Thus, we need to give conditions ensuring the
nonexplosion of $\{\xi_t,t\geq0\}$ [i.e., $P^\pi_x(\xi_t\in
S)\equiv1$]. To do so, we consider the following condition.
\renewcommand{\theassumm}{\Alph{assumm}}
\begin{assumm}\label{assumA}
There exist a continuous function $w\ge
1$ on $S$ and constants $\rho, b\ge0$ and a sequence of
nondecreasing subsets $\{S_k\}$ of $S$, such that:
\begin{longlist}[(3)]
\item[(1)] $\int_Sw(y)q(dy|x,a)\le\rho w(x)+b$ for
all $(x,a)\in K$;

\item[(2)] $ \inf_{x\notin S_k}w(x) \uparrow+\infty$ as
$k\to\infty$, with $\inf\varnothing:=\infty$;

\item[(3)] $ S_k\uparrow S$ and ${\sup_{a\in A(x), x\in
S_k}}|q(\{x\}|x,a)|<\infty$ for all $k\geq1$.
\end{longlist}
\end{assumm}
\begin{rem}\label{rem5.0-a}
We call Assumption \ref{assumA} a
nonexplosion condition for $\{\xi_t,t\geq0\}$. Obviously,
Assumption \ref{assumA} trivially holds when the transition rates are
\textit{bounded}; see \cite{havis,K71,K85,LP01,Pi98,P94,S99,Y77}, for
instance.
Assumption \ref{assumA} is similar to those in \cite{C04,G07b,G07a,GH09,GHP06}
for Markov policies and unbounded transition rates, and it can be
verified with
examples in
\cite{C04,G07b,G07a,GH09,GHP06} and those below.
\end{rem}

Under Assumption \ref{assumA}, we see (by Theorem \ref{t2} below) that
$\{\xi
_t,t\geq0\}$ is nonexplosive. Thus, for any fixed discount factor
$\alpha>0$ and an
initial distribution $\gamma\in{\mathcal{P}}(S)$, we define the
expected discounted criteria
%
%
\begin{eqnarray}\label{2.10}
V_\alpha(x,\pi,u)&:=&\int_0^\infty e^{-\alpha t}\int_A
E^\pi_x[u(\xi_{t-},a)\pi(da|e,t)]\,dt,\nonumber\\[-8pt]\\[-8pt]
V_\alpha(\pi,u)&:=&\int_SV_\alpha(x,\pi,u)\gamma(dx)\nonumber
\end{eqnarray}
for each $\pi\in\Pi, x\in S$ and a
measurable function $u$ on $K$, provided the integrals in
(\ref{2.10}) are well defined.

In particular, let
\[
V_r(x,\pi):=V_\alpha(x,\pi,r),\qquad V_r(\pi):=V_\alpha(\pi,r)
\]
and
\[
V_n(x,\pi):= V_\alpha(x,\pi,c_n),\qquad V_n(\pi):=
V_\alpha(\pi,c_n) \qquad\mbox{for } n=1, \ldots, N.
\]
[The finiteness of $V_r(\pi)$ and $V_n(\pi)$ will be ensured in
Theorem \ref{t3} below.]

Let
%
%
\begin{equation} \label{2.8}
U:=\{\pi| V_n(\pi)\leq d_n, n=1,\ldots, N\}\quad \mbox{and}\quad
V_r(U):=\sup_{\pi\in U} V_r(\pi)
\end{equation}
be the set of constrained policies and the constrained optimal reward
value, respectively.

In the following arguments, we assume that the set $U$ is \textit{not}
empty, and the discount factor $\alpha$ and
the initial distribution $\gamma$ as well as the numbers $d_n$ are
\textit{fixed}.

Then, the constrained optimality problem under consideration is as follows:
%
%
\begin{equation}\label{P}
\mbox{Maximize } V_r(\pi) \mbox{ over all } \pi\in U.
\end{equation}

\begin{defn}\label{defn2.2}
A policy $\pi^*\in U$ is said to be constrained optimal if $V_r(\pi
^*)=V_r(U)$. When $U=\Pi$, a constrained optimal policy is said to be
unconstrained optimal.
\end{defn}

The main goal of this paper is to give the conditions for the
existence and solvability of a constrained/unconstrained optimal policy.

\section{Main results}\label{section-3}

We state the main results of our work in this section. Their proofs are
presented later
in Section \ref{proof}. The main results are given
in three subsections.

\subsection{Conditions for nonexplosion and
finiteness}\label{section-5}

This subsection states the results on the nonexposition of
$\{\xi_t,t\geq0\}$ and finiteness of $V_n(x,\pi)$ and $V_n(\pi)$.

For the nonexposition of $\{\xi_t,t\geq0\}$, we have the following
fact.
\begin{theorem}\label{t2}
Suppose that Assumption \ref{assumA} holds. Then,
for each $\pi\in\Pi$, $x\in S$ and $t\ge0$:
\begin{longlist}[(a)]
\item[(a)] $P^\pi_x(T_\infty=\infty)=1$ and $P^\pi_x(\xi_t\in
S)=1$.\vspace*{2pt}

\item[(b)]
\[
E^\pi_x[w(\xi_t)]\leq
\cases{
\displaystyle e^{\rho t}w(x)+\frac{b}{\rho}(e^{\rho t}-1), &\quad if $\rho\not=0$,
\cr
w(x)+bt, &\quad if $\rho=0$.}
\]

\item[(c)] The analog of the forward Kolmogorov equation holds:
\[
P^\pi_x(\xi_t\in
D)=I_D(x)+E^\pi_x\biggl[\int_0^t\int_A\pi(da|e,s)q(D|\xi_{s-}(e),a)\,ds\biggr]
\]
for each $D\in{\mathcal{B}}(S)$ with $\sup_{x\in D}q^*(x)<\infty$.
\end{longlist}
\end{theorem}

The proof of Theorem \ref{t2} appears in Section \ref{proof}.
\begin{rem}
Theorem \ref{t2}(a) establishes the
nonexplosion of $\{\xi_t,t\geq0\}$ on the
probability space $(\Omega,\mathcal{F}, P_x^\pi)$ (for each policy
$\pi
\in\Pi$ and $x\in S$), and Theorem~\ref{t2} is an extension of the
corresponding results in \cite
{havis,K85,K71,LP01,Pi05,Pi98,P94,S99,Y77} for bounded transition rates and
in \cite{C04,G07b,G07a,GH09,GH03a,GH03c,GHP06,GR06,LMT96} for Markov
policies only. The process $\{\xi_t,t\geq0\}$ may \textit{not} be
Markovian because a policy $\pi$ can depend on state histories.
\end{rem}

Inspired by Theorem \ref{t2}, we introduce the following
condition.
\begin{assumm}\label{assumB}
Let $c_0(x,a):=-r(x,a)$ for $(x,a)\in
K$, and $w$ be as in Assumption \ref{assumA}.
\begin{longlist}[(3)]
\item[(1)] There exists a constant $M>0$ such that,
$|c_n(x,a)|\leq Mw(x)$ for every $(x,a)\in K$ and
$n=0,1,\ldots, N$.
\item[(2)] The discount factor $\alpha$ satisfies that
$\alpha>\rho$, with $\rho$ as in Assumption~\ref{assumA}.
\item[(3)] $\int_S
w(x)\gamma(dx)<\infty$.
\end{longlist}
Then the following fact establishes the finiteness of $V_n(x,\pi)$ and
$V_n(\pi)$.
\end{assumm}
\begin{theorem}\label{t3}
Suppose that Assumptions \ref{assumA} and \ref{assumB} hold.
Then, for each $\pi\in\Pi$ and $x\in S$:
\begin{longlist}[(a)]
\item[(a)] $E^\pi_x[|c_n(\xi_t,a)|\pi(da|e,t)]\leq
ME^\pi_x[w(\xi_t)]$ for all $t\geq0$ and $n=0,1,\ldots, N;$
\item[(b)] $|V_n(x,\pi)|\le M[\alpha w(x)+b]/[\alpha(\alpha-\rho)]$
and $|V_n(\pi)|\le MM^*_1$ for $n=0,1,\ldots, N$,
where $V_0(x,\pi):=V_\alpha(x,\pi,c_0), V_0(\pi):=V_\alpha(\pi,c_0),
M^*_1:=[\alpha\times\break\int_S
w(x)\gamma(dx)+b]/[\alpha(\alpha-\rho)]$.
\end{longlist}
\end{theorem}
\begin{pf}
Obviously, this theorem follows from Theorem
\ref{t2}(b) and (\ref{2.10}).
\end{pf}

\subsection{Existence of constrained optimal policies}

This subsection states the main results on the existence of
constrained optimal policies.

In order to show the existence of a
constrained optimal policy, as in \cite
{A99,AS91,OJ00,OJ03,OL96,KNH,Pi05,Pi97}, we introduce a key concept of
an occupation
measure of a policy.
\begin{defn}\label{De-4.1}
Fix policies $\pi,\pi_1,\pi
_2\in\Pi$.
\begin{longlist}[(iii)]
\item[(i)] The occupation measure of $\pi$ is a probability
measure $\eta^\pi$ on $S\times A$ concentred on $K$, which is
defined by
%
%
\begin{eqnarray}\label{4.1}
\eta^\pi(D\times\Gamma):=\alpha\int_0^\infty e^{-\alpha
t}E^\pi_\gamma\bigl[I_{\{\xi_t\in
D\}}(e)\pi(\Gamma|e,t)\bigr]\,dt\nonumber\\[-8pt]\\[-8pt]
&&\eqntext{\mbox{with } D\in{\mathcal{B}}(S),
\Gamma\in{\mathcal{B}}(A).}
\end{eqnarray}
(Obviously, $\eta^\pi$ concentrates on $K$ and depends on $\pi,
\alpha$ and $\gamma$. However, we impress $\gamma$ and $\alpha$ in
the occupation measure for simplicity.)
\item[(ii)] Two policies
$\pi^1$ and $\pi^2$ are called equivalent if
$\eta^{\pi^1}=\eta^{\pi^2}$.
\item[(iii)] We denote by $\hat\eta$
the \textit{marginal} (or \textit{projection}) on $S$ of
a probability measure $\eta$ on $S\times A$, and by $\phi^\eta(\in
\Pi
_s)$ the randomized stationary policy
(depending on $\eta$), which is determined by the following
decomposition of
$\eta$:
%
%
\begin{equation}\label{4.2}
\eta(dx,da)=\hat\eta(dx)\phi^\eta(da|x).
\end{equation}
\end{longlist}
\end{defn}

Thus, by (\ref{4.1}) and (\ref{2.10}), we have
$V_\alpha(x,\pi,u)=\frac{1}{\alpha}\int_{S\times
A}u(x,a)\eta^\pi(dx,da)$, and we can rewrite (\ref{P}) as an
\textit{equivalent} optimality problem:
%
%
\begin{eqnarray}\label{4.2b}
&&\mbox{Maximize } \frac{1}{\alpha}\int_{K}
r(x,a)\eta(dx,da)\nonumber\\[-8pt]\\[-8pt]
&&\qquad\mbox{over } \eta\in\biggl\{\eta^\pi\dvtx\int_{K} c_n(x,a)\eta^\pi
(dx,da)\leq
\alpha d_n, 1\leq n\leq N\biggr\}.\nonumber
\end{eqnarray}
To solve problem (\ref{4.2b}), we need to seek a certain
compactness structure on the set of all occupation measures. To do
so, we require to characterize an occupation measure, and we have the following
fact.
\begin{theorem}\label{Th-4.1}
Under Assumption \ref{assumA}, the
following assertions hold.
\begin{longlist}[(a)]
\item[(a)] The occupation measure $\eta^\pi$ (for each fixed $\pi
\in
\Pi$) satisfies
the following equation:
\begin{eqnarray}
\alpha\hat\eta^\pi(D)=\alpha\gamma(D)+\int_{S\times A}
q(D|x,a)\eta^\pi(dx,da) \nonumber\\
&&\eqntext{\displaystyle \forall D\in{\mathcal{B}}(S)
\mbox{ with } \sup_{x\in D} q^*(x)<\infty.}
\end{eqnarray}
\item[(b)] Conversely, if a probability measure $\eta$ on
$S\times A$ (concentrated on $K$) satisfies
\begin{eqnarray}
\alpha\hat\eta(D)=\alpha\gamma(D)+\int_{S\times
A}q(D|x,a)\eta(dx,da) \nonumber\\
&&\eqntext{\displaystyle \forall D\in{\mathcal{B}}(S) \mbox{ with }
\sup_{x\in D} q^*(x)<\infty}
\end{eqnarray}
and $\int_S |q(\{x\}|x,\phi^\eta)|\hat\eta(dx)<\infty$, then
$\eta^{\phi
^\eta}=\eta$, where $\phi^\eta$ is as in (\ref{4.2}).
\item[(c)] If, in addition, Assumptions \ref{assumB}(2) and \ref
{assumB}(3) are
satisfied, and $q^*(x)\leq Lw(x)$ for all $x\in S$, with some
constant $L>0$, then $\phi^{\eta^\phi}=\phi$ for all $\phi\in
\Pi_s$.
\end{longlist}
\end{theorem}

The proof of Theorem \ref{Th-4.1} appears in
Section \ref{proof}.
\begin{rem}\label{rem5.1-r}
Theorems \ref{Th-4.1}(a) and \ref{Th-4.1}(b)
are proved in
\cite{Pi05} for continuous-time MDPs with uniformly bounded
transition rates and in \cite{A99,AS91,OL96} for discrete-time
MDPs.
\end{rem}

To give a certain convergence of occupation measures, we introduce
some notation.

For any real-valued continuous function $\bar w\geq1$ on $S$, let
\[
{\mathcal{P}}_{\bar w}(S\times A):=\biggl\{\eta\in{\mathcal{P}}(S\times A)\Big|
\int_S\bar w(x)\hat\eta(dx)<\infty\biggr\}.
\]
Then we define two maps, $T_{\bar w}$ and $T_{\bar w}'$, as follows:
\[
T_{\bar w}:\quad{\mathcal{P}}_{\bar w}(S\times A) \longrightarrow
{\mathcal{P}}(S\times A),\qquad \eta\mapsto T_{\bar w}(\eta),
\]
where $T_{\bar w}(\eta)$ is given by
%
\begin{equation}
\label{4.13}
T_{\bar w}(\eta)(D\times\Gamma):=\frac{\int_D\bar
w(x)\eta(dx,\Gamma)}{\int_S\bar w(x)\hat\eta(dx)} \qquad\forall
D\in{\mathcal{B}}(S) \mbox{ and } \Gamma\in{\mathcal{B}}(A);
\end{equation}
\[
T_{\bar w}':\quad{\mathcal{P}}(S\times A) \longrightarrow
{\mathcal{P}}_{\bar w}(S\times A),\qquad \mu\mapsto T_{\bar
w}'(\mu),
\]
where $T_{\bar w}'(\mu)$ is given by
%
\begin{equation}\qquad
\label{4.14}
T_{\bar w}'(\mu)(D\times\Gamma):=\frac{\int_D ({1}/{\bar
w(x)})\mu
(dx,\Gamma)}{\int_S({1}/{\bar w(x)})\hat\mu(dx)}
\qquad\forall D\in{\mathcal{B}}(S) \mbox{ and }
\Gamma\in{\mathcal{B}}(A).
\end{equation}
[Since $1\leq\bar w<\infty$ on $S$, we have $0<\int_S\frac{1}{\bar
w(x)}\mu(dx)\leq1$ for any $\mu\in{\mathcal{P}}(S)$, and thus the
maps $T_{\bar w}$ and $T_{\bar w}'$ are well defined.]
\begin{defn}\label{De-4.2}
The $\bar w$-weak topology on
${\mathcal{P}}_{\bar w}(S\times A)$ is defined by the $\bar w$-weak
convergence as follows: a
sequence $\{\eta_k, k\geq1\}\subseteq{\mathcal{P}}_{\bar w}(S\times
A)$ is called to $\bar w$-converge weakly to $\eta\in
{\mathcal{P}}_{\bar w}(S\times A)$ (and written as $\eta_k \stackrel
{\bar w}{\longrightarrow}\eta$) if
\[
\lim_{k\to\infty}
\int_{S\times A} u(x,a)\eta_k(dx,da)=\int_{S\times A} u(x,a)\eta(dx,da)
\]
for each continuous function $u(x,a)$ on
$S\times A$ such that $|u(x,a)|\leq L_u \bar w(x)$ for all $(x,a)\in
S\times A$, with some nonnegative constant $L_u$ depending on $u$.
\end{defn}

Obviously, $\eta_k \stackrel{\bar
w}{\longrightarrow}\eta$ implies $\eta_k \stackrel
{1}{\longrightarrow}
\eta$ (the standard weak convergence of probability measures).
The following lemma establishes the relationship between $\bar w$- and
standard weak convergence.
\begin{lem}\label{Th-4.5}
For any given real-valued continuous
function $\bar w\geq1$ on $S$, let
$\{\eta_k, k=0,1,\ldots\} \subset{\mathcal{P}}_{\bar w}(S\times A)$
and $\{\mu_k, k=0,1,\ldots\} \subset{\mathcal{P}}(S\times A)$. Then:
\begin{longlist}[(a)]
\item[(a)] $T_{\bar w}(\eta)\in{\mathcal{P}}(S\times A)$ for all $
\eta\in{\mathcal{P}}_{\bar w}(S\times A)$ and $T_{\bar
w}'(\mu)\in{\mathcal{P}}_{\bar w}(S\times A)$ for all $\mu\in
{\mathcal{P}}(S\times A)$;
\item[(b)] $T_{\bar w}'(T_{\bar w}(\eta))=\eta$ for all
$\eta\in{\mathcal{P}}_{\bar w}(S\times A)$ and
$T_{\bar w}(T_{\bar w}'(\mu))=\mu$ for all $\mu\in{\mathcal
{P}}(S\times A)$;
\item[(c)] $\eta_k \stackrel{\bar w}{\longrightarrow}\eta_0$ if
and only if $ T_{\bar
w}(\eta_k)\stackrel{1}{\longrightarrow}T_{\bar w}(\eta_0)$;
\item[(d)] $\mu_k \stackrel{1}{\longrightarrow}\mu_0$ if and only
if $
T_{\bar w}'(\mu_k)\stackrel{\bar
w}{\longrightarrow}T_{\bar w}'(\mu_0)$.
\end{longlist}
\end{lem}

The proof of Lemma \ref{Th-4.5} appears in
Section \ref{proof}.

To further analyze the properties of occupation
measures, we let
%
%
\begin{eqnarray}
\label{4.9}
{\mathcal{M}}_o&:=&\biggl\{\eta^\pi\Big| \int_Sw(x)\hat\eta^\pi
(dx)<\infty,
\pi\in\Pi\biggr\} \subseteq{\mathcal{P}}_w(K)\nonumber\\[-8pt]\\[-8pt]
&&\eqntext{\mbox{(with } w \mbox{ as
in Assumption \ref{assumA}),}} \\
\label{D}
{\mathcal{M}}_o^c&:=&\biggl\{\eta\in{\mathcal{M}}_o\Big| \int_{S\times A}
c_n(x,a)\eta(dx,da) \leq\alpha d_n, n=1, \ldots, N\biggr\}.
\end{eqnarray}
\begin{lem}\label{Th-4.4}
Suppose that Assumptions \ref{assumA}, \ref{assumB}(2)
and \ref{assumB}(3) hold. If, in addition, $q^*(x)\leq Lw(x)$ for all
$x\in S$,
with some constant $L>0$, then the following assertions hold:
\begin{longlist}[(a)]
\item[(a)] ${\mathcal{M}}_o$ and ${\mathcal{M}}_o^c$ are convex.

\item[(b)] If, in addition, $\int_Sg(y)q(dy|x,a)$ is continuous on
$K$ for each fixed $g\in
C_b(S)$, then ${\mathcal{M}}_o$ is closed (with respect to the
$w$-weak topology).
\end{longlist}
\end{lem}

The proof of Lemma \ref{Th-4.4} appears
in Section \ref{proof}.

For the solvability of
(\ref{4.2b}), by Lemmas \ref{Th-4.5} and \ref{Th-4.4}, we
introduce the following condition.
\begin{assumm}\label{assumC}
Let $w$ be as in Assumption \ref{assumA}.
\begin{longlist}[(3)]
\item[(1)] The functions $c_n(x,a)$ and $\int_Sg(y)q(dy|x,a)$
are continuous on $K$ [for each fixed $g\in
C_b(S)$ and $0\leq n\leq N$].
\item[(2)] There exist a measurable function $w'\geq1$ on $S$
and a nondecreasing sequence of compact sets
$K_m\uparrow K$, such that $\lim_{m\to\infty} \inf_{(x,a)\notin
K_m}\frac{w(x)}{w'(x)}=\infty$. 
\item[(3)] There exist a constant $L>0$ such that $q^*(x)\leq
Lw(x)$ for all $x\in S$.
\end{longlist}
\end{assumm}
\begin{rem}\label{rem5.1-new}
Assumption \ref{assumC}(2) is slightly different
from the compactness condition in \cite{OJ00,OJ03,OL99,OL96,KNH} for
discrete-time MDPs and \cite{G07b,GR06} for
continuous-time MDPs.
\end{rem}

We now state our second main result on the existence of a
constrained optimal policy.
\begin{theorem}\label{Th-4.6}
Suppose that Assumptions \ref{assumA}, \ref{assumB}
and \ref{assumC} hold. Then:
\begin{longlist}[(a)]
\item[(a)] ${\mathcal{M}}_o$ and ${\mathcal{M}}_o^c$ are metrizable and
compact (with respect to the $w'$-weak topology), that is,\vspace*{-1pt} for any sequence
$\{\eta_k,k\geq1\}$ in ${\mathcal{M}}_o$ (or ${\mathcal{M}}_o^c$), there
exists a subsequence $\{\eta_{k_m},m\geq1\}$ and
$\eta_0\in{\mathcal{M}}_o$ (or ${\mathcal{M}}_o^c$) such that such that
$\eta_{k_m} \stackrel{w'}{\longrightarrow} \eta_0$ as $m\to\infty$.
\item[(b)] There exists a constrained optimal policy.
\end{longlist}
\end{theorem}

The proof of Theorem \ref{Th-4.6} appears
in Section \ref{proof}.
\begin{rem}\label{rem5.1}
Theorem \ref{Th-4.6}(b) shows the
existence of a constrained optimal policy. It
should be noted that the conditions for Theorem \ref{Th-4.6}(b) are
weaker than those in \cite{G07b,GH09,GH03a,GH03c,P94} for the class of
all Markov policies. This is because some assumptions such as the
nonnegativity of costs in
\cite{GH03c} and the absolute integrability condition in \cite
{G07b,GH09,GH03c} are not required here.
\end{rem}

\subsection{Solvability of constrained optimal policies}

This subsection states the results on the solvability of constrained
optimal policies.

First, by (\ref{4.2b}) we see that the original constrained
optimality problem (\ref{P}) is \textit{equivalent} to the
following constrained minimization problem:
%
%
\begin{equation} \label{5.1}
\mbox{Minimize } V_0(\pi) \mbox{ over } \pi\in\{\pi| V_n(\pi
)\leq d_n, n=1,\ldots, N\}.
\end{equation}

By (\ref{2.10}) and (\ref{4.1}), the problem
(\ref{5.1}) can be rewritten into the following form:
\[
\cases{
\displaystyle\inf_{\eta\in\{\eta^\pi| \pi\in\Pi\}}
\frac
{1}{\alpha} \int_{S\times A}
c_0(x,a)\eta(dx,da), \vspace*{2pt}\cr
\displaystyle\mbox{subject to } \int_{S\times A} c_n(x,a)\eta(dx,da)
\leq
\alpha d_n, &\quad$n=1, \ldots, N$,}
\]
which (by Theorem \ref{Th-4.1}) is equivalent to the following linear
program (LP):
%
%
\begin{equation}
\label{LP}
\mbox{LP}:\quad \inf_{\eta} \int_{S\times A} \frac{1}{\alpha}
c_0(x,a)\eta(dx,da)
\end{equation}
subject to
%
{\setcounter{equation}{8}
\renewcommand{\theequation}{\arabic{section}.\arabic{equation}$'$}
\begin{equation}
\label{LP'}
\cases{\displaystyle \int_{S\times A} c_n(x,a)\eta(dx,da) \leq
\alpha d_n, \qquad n=1,
\ldots, N,\vspace*{2pt}\cr
\displaystyle \alpha\hat\eta(D)=\alpha\gamma(D)+\int_{S\times A} q(D|x,a)\eta
(dx,da),\vspace*{2pt}\cr
\qquad\mbox{for all } D\in{\mathcal{B}}(S) \mbox{
with $\displaystyle \sup_{x\in D} q^*(x)<\infty$}, \vspace*{2pt}\cr
\displaystyle \int_{S}w(x)\hat\eta(dx)<\infty, \qquad\hspace*{48.7pt} \eta\in{\mathcal{P}}(K).}
\end{equation}
}

\vspace*{-4pt}

\noindent Obviously, (\ref{LP}) is a linear program over the set of probability
measures $\eta\in{\mathcal{P}}(K)$ satisfying
(\ref{LP'}). We call (\ref{LP}) the primal \textit{linear programming
formulation} of (\ref{P}).

Thus, we obtain the following result on the solvability of
constrained optimal policies.
\begin{theorem}\label{Th-5.1}
Under Assumptions \ref{assumA}, \ref{assumB} and \ref{assumC}(3),
the following assertions hold.
\begin{longlist}[(a)]
\item[(a)] If there exists a feasible solution to LP (\ref{LP}), then
the set $U$ of
constrained policies is nonempty. Conversely, if $U$ is nonempty, then
there exists a feasible solution to LP
(\ref{LP}).
\item[(b)] If there exists an optimal solution $\eta^*$ to LP (\ref
{LP}), then the randomized stationary
policy\vspace*{1pt} $\phi^{\eta^*}$ is constrained optimal. Conversely, if $\pi^*$
is constrained optimal, then
$\eta^{\pi^*}$ is an optimal solution to LP (\ref{LP}).
\item[(c)] If, in addition, $U\not=\varnothing$ and Assumptions \ref
{assumC}(1) and
\ref{assumC}(2) are satisfied, then an optimal solution $\eta^*$ to LP
(\ref{LP})
exists, and the policy $\phi^{\eta^*}$ is constrained optimal.
\end{longlist}
\end{theorem}

The proof of Theorem \ref{Th-5.1} appears
in Section \ref{proof}.

In particular, when $S$ and $A(x)$ are finite, then LP
(\ref{LP}) is the form of
%
%
\begin{eqnarray}\label{LP-f}\qquad
&& \mbox{minimize } \sum_{x\in S}\sum_{a\in A(x)} \frac{1}{\alpha}
c_0(x,a)\eta(x,a) \nonumber\\[-8pt]\\[-8pt]
&&\qquad\mbox{subject to } \cases{
\displaystyle \sum_{x\in S}\sum_{a\in A(x)} c_1(x,a)\eta(x,a)\leq\alpha d_1, \cr
\ \ \ \ \ \vdots  \ \ \ \ \ \ \ \ \ \  \ \ \ \ \ \vdots  \ \ \ \ \ \ \ \ \ \  \ \ \ \ \ \vdots \cr
\displaystyle \sum_{x\in S}\sum_{a\in A(x)} c_n(x,a)\eta(x,a) \leq\alpha d_N,
\cr
\displaystyle \alpha\sum_{a\in A(x)}\eta(x,a)=\alpha\gamma(x)+\sum_{y\in
S}\sum_{a\in A(y)} q(x|y,a)\eta(y,a),\vspace*{2pt}\cr
\qquad\forall x\in S, \eta(x,a)\geq0, x\in S, a\in A(x),}\nonumber
\end{eqnarray}
which is a LP and can be solved by many methods such as the well-known
simplex method.

To state the structure of constrained optimal policies, we need to
recall some concepts. We say that under $\phi\in\Pi_s$, there are
$m(x,\phi)$ randomizations at $x\in S$ if there are $m(x,\phi)+1$
actions $a\in A(x)$ for which $\phi(a|x)>0$. When $S$ and $A(x)$
are finite, we call $\#(\phi):=\sum_{x\in S}m(x,\phi)$ the number
of randomizations under $\phi$.

Thus, following Theorem 3.8 in \cite{A99} and Theorem \ref{Th-5.1}
above, we have the following fact.
\begin{cor}\label{Co-1}
Suppose that $S$ and $A(x)$ are
finite. Let $\eta^*$ be an optimal basic solution to LP (\ref{LP-f}).
Then, the policy $\phi^{\eta^*}$ is constrained optimal, where
$\phi^{\eta^*}$ is given by
%
%
\begin{equation}\label{mr}
\phi^{\eta^*}(a|x)=\cases{
\displaystyle \frac{\eta^*(x,a)}{\hat\eta^*(x)}, &\quad when
$\displaystyle \hat\eta^*(x):=\sum_{a\in A(x)}\eta^*(x,a)>0$\cr
&\quad and $a\in A(x)$,\cr
I_{\{a(x)\}}(a), &\quad when $\hat\eta^*(x)=0$
\mbox{ and } $a \in A(x)$,}
\end{equation}
for all $x\in S$, $a(x)\in A(x)$ is chosen arbitrarily. Further, $\#
(\phi^{\eta^*})\leq N$.
\end{cor}

Corollary \ref{Co-1} provides the structure of a constrained optimal
policy for finite $S$ and $A(x)$, and it is proven for the case of
denumerable states and a single constraint in \cite{GH03c,ZG08}. For
a more general case of Polish spaces, we have the following facts,
in which the first one (i.e., Theorem \ref{Th-3.2}) establishes the
relationship between stationary policies in $F$ and extreme points
in ${\mathcal{M}}_o$, and the second one (i.e., Theorem \ref{Th-5.2})
shows a structure of a constrained optimal policy.
\begin{theorem}\label{Th-3.2}
Suppose that Assumptions \ref{assumA},
\ref{assumB}(2), \ref{assumB}(3) and \ref{assumC}(3) hold. Then:
\begin{longlist}[(a)]
\item[(a)] $\eta^f$ is an extreme point in ${\mathcal{M}}_o$ for each
$f\in F$.
\item[(b)] If, for each $\phi\in\Pi_s$ and $D\in{\mathcal{B}}(S)$ with
$\hat\eta^\phi(D)>0$, there exists state
$x \in D$ (depending on $D$ and $\phi$) such that $\hat\eta^\phi(\{
x\}
)>0$, then $\eta$ is an extreme point in ${\mathcal{M}}_o$
\textit{if} and \textit{only if} there exists a policy $f \in F$ such that
$\eta=\eta^f$.
\end{longlist}
[\textup{The condition in Theorem \ref{Th-3.2}(b) is satisfied when
$S$ is denumerable.}]
\end{theorem}

The proof Theorem \ref{Th-3.2} appears in
Section \ref{proof}.
\begin{theorem}\label{Th-5.2}
Suppose that Assumptions \ref{assumA}, \ref{assumB},
\ref{assumC} and the conditions for Theorem \ref{Th-3.2}\textup{(b)} are
satisfied. Then,
there exists a
constrained optimal policy $\pi^*\in\Pi_s$, which is a mixture of
$(N+1)$ stationary policies, that is, there exists $(N+1)$ numbers
$p_n\geq0$ and policies $f_n\in F$ $(1\leq n\leq N+1)$ such that $\pi
^*=\phi^{(p_1 \eta^{f_1}+\cdots+p_{N+1} \eta^{f_{N+1}})}$ and
$p_1+\cdots+p_{N+1}=1$.
\end{theorem}

The proof of Theorem \ref{Th-5.2} appears in
Section \ref{proof}.
\begin{rem}\label{rem5.1-rr}
The arguments of Theorems
\ref{Th-4.6}, \ref{Th-5.1}, \ref{Th-3.2} and \ref{Th-5.2} do not
depend on the data in model (\ref{2.1}), but they are based on
Theorem \ref{Th-4.1}. Thus, the discrete-time versions of Theorems
\ref{Th-4.6}, \ref{Th-5.1}, \ref{Th-3.2} and \ref{Th-5.2} are
still true because Theorem \ref{Th-4.1} is established in
\cite{A99,AS91,OL96} for discrete-time MDPs.
\end{rem}

\section{Examples}\label{section-4}

In this section, we illustrate our conditions and main results
with examples.
\begin{exm}\label{exm5.1} Let $S:=(-\infty,\infty)$,
$A(x):=[\beta
_0, \beta(|x|+1)]$ for each
$x\in S$ with some constants $0<\beta_0< \beta$. Suppose that the
reward $r(x,a)$ and costs $c_n(x,a)$ $(1\leq n\leq N)$ are given.
We consider the transition
rates $q(\cdot|x,a)$ given by
%
%
\begin{eqnarray}\label{e-1}
q(D|x,a):=(|x|+1)\biggl[\int_{D-\{x\}}f(y|x,a)\,dy-\delta_x(D)
\biggr] \nonumber\\[-8pt]\\[-8pt]
&&\eqntext{\mbox{for }
(x,a)\in K, D\in{\mathcal{B}}(S),}
\end{eqnarray}
where $f(y|x,a):= \frac{1}{\sqrt{2\pi a}}e^{-{(y-x)^2}/({2a})}$
is the density function of Gaussian distribution $N(x,a)$.
\end{exm}

We now aim to find conditions that ensure the existence of constrained
optimal policies for Example \ref{exm5.1}.
To do so, we need the following hypotheses.
\begin{assumm}\label{assumD}
Let $\alpha, \gamma, d_n$ and $U$ $(\mbox{$\not=$}\varnothing)$
be as in (\ref{2.8}).
\begin{longlist}[(2)]
\item[(1)] $\alpha>6\beta$ and $\int_S
x^4\gamma(dx)<\infty$ (hence, there exists a constant $\rho$ such
that $6\beta<\rho<\alpha$);
\item[(2)] $c_n(x,a)$ ($0 \leq n\leq N$) are continuous on $K$ and
$|c_n(x,a)|\leq L'(x^2+1)$ for all $(x,a)\in K$,
with some constant $L'>0$, where $c_0(x,a):=-r(x,a)$.
\end{longlist}
\end{assumm}

Then, we have the following result.
\begin{prop}\label{prop-exm01} Under Assumption \ref{assumD}, Example
\ref{exm5.1}
satisfies Assumptions~\ref{assumA}, \ref{assumB} and
\ref{assumC}. Therefore (by Theorem \ref{Th-4.6}), there exists a constrained
optimal policy for Example \ref{exm5.1}.
\end{prop}
\begin{pf}
For each $m\geq1$ and $x\in S$, let
%
%
\begin{eqnarray}\label{e-3}
S_m&:=&[-m,m],\qquad K_m:=\{(x,a)| x\in S_m, a\in
A(x)\},\nonumber\\[-8pt]\\[-8pt]
w'(x)&:=&x^2+1,\qquad w(x):=x^4+1.\nonumber
\end{eqnarray}
To verify Assumption \ref{assumA}, it suffices to verify Assumption \ref
{assumA}(1)
because Assumptions \ref{assumA}(2) and \ref{assumA}(3) follow from
(\ref{e-3}) and
(\ref{e-1}). Indeed, by (\ref{e-1}) and a straightforward
calculation, we have
%
%
\begin{eqnarray}\label{e-4}
\int_S w(y)q(dy|x,a)
&=&6(x^2 a+3a^2)(|x|+1) \nonumber\\[-8pt]\\[-8pt]
&\le&\beta w(x)+b\qquad \mbox{for some constant }b>0,\nonumber
\end{eqnarray}
which implies Assumption \ref{assumA}(1).

Obviously, Assumption \ref{assumB} follows from (\ref{e-4}) and
Assumptions \ref{assumD}(1)
and \ref{assumD}(2).

To verify Assumption \ref{assumC}, for any $g\in C_b(S)$, by (\ref
{e-1}) we
have the following:
\[
\int_S g(y)q(dy|x,a)=(|x|+1)\biggl[\int_{-\infty}^\infty g(y)\frac
{1}{\sqrt{2\pi a}}e^{-{(y-x)^2}/({2a})}\,dy-g(x)\biggr],
\]
which, together with the dominated convergence theorem, implies
Assumption~\ref{assumC}(1). Therefore, Assumption \ref{assumC} holds
because Assumptions
\ref{assumC}(2) and \ref{assumC}(3) follow from (\ref{e-1}) and (\ref{e-3}).

Using Example \ref{exm5.1}, we present computable examples for
unconstrained optimal policies.
\begin{exm}\label{exm5.2} With the same data as in Example \ref
{exm5.1}, we further suppose that $r(x,a)$ in
Example \ref{exm5.1} is given by
%
%
\begin{equation}\label{e-2}
r(x,a):= px^2-\delta a^2\qquad \mbox{for } (x,a)\in K,
\end{equation}
where $p,\delta>0$ are
fixed constants.
\end{exm}
\begin{assumm}\label{assumE}
Let $\beta_0$ and $\beta$ be as in
Example \ref{exm5.1}, and $L'$ as in Assumption~\ref{assumD}(2).
\begin{longlist}[(2)]
\item[(1)] $d_n\geq L'[\alpha\int_Sx^4\gamma(dx)+\alpha+b]/[\alpha
(\alpha-\beta)]$ for all $1\leq n\leq
N)$, with $ b:=\beta(\frac{\rho+2\beta}{\rho-\beta}+2)^2$;
\item[(2)] $2\alpha\beta_0-\beta_0^2\le\frac{p}{\delta}\le\min
\{
\alpha^2,2\alpha\beta-\beta^2\}
$, with $p, \delta$ as in (\ref{e-2}).
\end{longlist}
\end{assumm}
\begin{prop}\label{prop-exm01-a} Suppose that Assumptions \ref{assumD}
and \ref{assumE}
hold. Then:
\begin{longlist}[(a)]
\item[(a)] Example \ref{exm5.2} satisfies Assumptions \ref{assumA}, \ref
{assumB} and \ref{assumC}.
Moreover, $V_r(U)=\int_Su(x)\gamma(dx)$, where
\begin{eqnarray*}
u(x)&=& \bigl(2\delta\alpha- 2\sqrt{\delta^2\alpha^2-p\delta}\bigr) x^2+
\biggl(4\delta
\alpha- 4\sqrt{\delta^2\alpha^2-p\delta}
-\frac{2p}{\alpha}\biggr) |x|\\
&&{} +2\delta\alpha- 2\sqrt{\delta^2\alpha
^2-p\delta}- \frac{p}{\alpha}.
\end{eqnarray*}
\item[(b)] The stationary policy $f^*$ is unconstrained optimal for
Example \ref{exm5.2}, where
\[
f^*(x):= \Biggl(\alpha-\sqrt{\alpha^2-\frac{p}{\delta}}\Biggr)
(|x|+1)\qquad
\forall x\in S.
\]
\end{longlist}
\end{prop}
\begin{pf}
Note that Assumptions \ref{assumE}(1) and \ref{assumD} imply that $U=\Pi$
(by Theorem~\ref{t3}), and so the problem (\ref{P}) becomes an
unconstrained optimality problem. Thus, as in Proposition
\ref{prop-exm01}, under Assumptions \ref{assumD} and \ref{assumE}, we
see that all
assumptions in Theorem 3.3 in \cite{G07b} are satisfied. Hence,
Theorem 3.3 in \cite{G07b} ensures the existence of a function
$u$ in $B_w(S)$ such that, for each $x\in S$ and $\pi\in\Pi$,
%
%
\begin{equation}\label{e-06}\qquad
\alpha u(x)=\sup_{a\in A(x)} \biggl\{r(x,a) + \int_{S}
u(y)q(dy|x,a)
\biggr\}\quad \mbox{and}\quad u(x)\geq V_r(x,\pi).
\end{equation}
To obtain the analytic expression of $u$, we assume for a moment that
%
%
\begin{equation}\label{V}
\qquad u(x):= l_2 x^2+l_1x+l_0 \qquad\mbox{for } x\in S\mbox{,  with some
constants }
l_1,l_2,l_2.
\end{equation}
Then, using (\ref{e-1}), (\ref{e-2}) and (\ref{e-06}), by a straightforward
calculation we have
%
%
\begin{equation}\label{e-07}\quad
\alpha(l_2x^2+l_1 x + l_0)=\sup_{a\in A(x)} \biggl\{ px^2 - \delta
\biggl(a-\frac{ l_2(|x|+1)}{2\delta}\biggr)^2+
\frac{l_2^2(|x|+1)^2}{4\delta}\biggr\},\hspace*{-24pt}
\end{equation}
which implies that $f^*(x):=\frac{ l_2 (|x|+1)}{2\delta}$ attains the
maximum of the right-hand side of (\ref{e-07}).
Therefore, by Theorem 3.3 in
\cite{G07b}, we have
%
%
\begin{eqnarray}\label{07}
V_r(x,f^*)=u(x)\quad \mbox{and}\quad \alpha(l_2x^2+l_1 x + l_0)
=px^2 + \frac{l_2^2(|x|+1)^2}{4\delta} \nonumber\\[-8pt]\\[-8pt]
&&\eqntext{\forall x\in S.}
\end{eqnarray}
Comparing with the coefficients of both sides in (\ref{07}), we obtain
%
%
\begin{equation}\label{08}
\alpha l_2= p+ \frac{l_2^2}{4\delta},\qquad
\alpha l_1= \cases{
\dfrac{l_2^2}{2\delta}, &\quad if $x\ge0$, \vspace*{2pt}\cr
-\dfrac{l_2^2}{2\delta}, &\quad otherwise,}\qquad
\alpha l_0= \frac{l_2^2}{4\delta}.
\end{equation}
Under Assumption \ref{assumE}, solving the system of equations (\ref
{08}) gives
\begin{eqnarray*}
l_2&=& 2\delta\alpha- 2\sqrt{\delta^2\alpha^2-p\delta},\qquad l_0=
2\delta\alpha- 2\sqrt{\delta^2\alpha^2-p\delta}-
\frac{p}{\alpha},\\
l_1&=& \cases{
\displaystyle 4\delta\alpha- 4\sqrt{\delta^2\alpha^2-p\delta}-\frac{2p}{\alpha},
&\quad
if $x\ge0$, \vspace*{2pt}\cr
\displaystyle -\biggl( 4\delta\alpha- 4\sqrt{\delta^2\alpha^2-p\delta}-\frac
{2p}{\alpha}\biggr), &\quad otherwise,}
\end{eqnarray*}
which, together with (\ref{V}) and (\ref{07}), yields
\begin{eqnarray*}
u(x)&=& \bigl(2\delta\alpha- 2\sqrt{\delta^2\alpha^2-p\delta}\bigr) x^2+
\biggl(4\delta
\alpha- 4\sqrt{\delta^2\alpha^2-p\delta}
-\frac{2p}{\alpha}\biggr) |x|\\
&&{} +2\delta\alpha- 2\sqrt{\delta^2\alpha
^2-p\delta}-
\frac{p}{\alpha},
\\
f^*(x)&=&\Biggl(\alpha-\sqrt{\alpha^2-\frac{p}{\delta}}\Biggr) (|x|+1)
\in
A(x)\quad \mbox{and}\quad V_r(x,f^*)=u(x)\qquad \forall x\in S.
\end{eqnarray*}
This, together with (\ref{e-06}) and (\ref{2.10}), completes the proof
of this
proposition.
\end{pf}
\begin{exm}\label{exm5.3} Let $S:=(-\infty,\infty)$, $A(x):=[0,
\beta(|x|+1)]$ for each
$x\in S$ with some constant $\beta>0$, and the reward $r(x,a)$ and
transition
rates $q(\cdot|x,a)$ are defined as follows: for each $(x,a)\in K$ and
$D\in
{\mathcal{B}}(S)$,
\begin{eqnarray*}
q(D|x,a)&:=&(\beta|x|+a)\biggl[\int_{D-\{x\}}\frac{1}{\sqrt{2\pi
(\beta(|x|+1)-a+1)}}\\
&&\hspace*{81.5pt}{}\times e^{-{(y-x)^2}/({2(\beta
(|x|+1)-a+1)})}\,dy-\delta
_x(D)\biggr].
\\
r(x,a)&:=&p|x|a- \delta a^2 \qquad\mbox{for } (x,a)\in K\mbox{,  with }
p, \delta>0.
\end{eqnarray*}
\end{exm}
\renewcommand{\theassumm}{E}
\begin{assumm}\label{assumEE}
$\alpha>\beta^2$; $\int_S
x^2\gamma(dx)<\infty$; and $\beta\ge\max\{1, \frac{p}{2\delta}\}$.
\end{assumm}

Then as the arguments for Example \ref{exm5.2} in Proposition \ref
{prop-exm01-a}, we have the
following results.
\begin{prop}\label{prop-exm02}
Under Assumption \ref{assumEE}, Example \ref{exm5.3} satisfies
Assumptions~\ref{assumA}, \ref{assumB}
and \ref{assumC}. Moreover,
if, in addition, $U=\Pi$, then $V_r(U)=\int_S u(x)\gamma(dx)$, where
\begin{eqnarray*}
u(x)&=&\frac{1}{2}\delta\bigl(\sqrt{\kappa+1}-1\bigr)x^2\\
& &{} +
\frac{1}{2\alpha\kappa}\bigl[p\bigl(\sqrt{\kappa+1}-1\bigr)+\kappa\delta
\beta
\bigr](\beta+1)\bigl(\sqrt{\kappa+1}-1\bigr)|x|\\
&&{} +\frac{1}{8\alpha\kappa}\delta(\beta+1)^2\bigl(\sqrt{\kappa+1}-1\bigr)^3
\end{eqnarray*}
with $\kappa:=\frac{p^2}{\delta^2(\alpha-\beta^2)}>0$, and the
following stationary policy $f^*$ is unconstrained optimal:
\[
f^*(x):=\frac{p(\sqrt{\kappa+1}-1)}{\delta\kappa}|x|+ \frac
{1}{2\kappa
}(\beta+1)\bigl(\sqrt{\kappa+1}-1\bigr)^2\qquad \forall x\in S.
\]
\end{prop}
\begin{pf}
The proof of Proposition \ref{prop-exm02} is similar
to
that of Proposition \ref{prop-exm01}, and thus the details are
omitted here.
\end{pf}
\noqed
\end{pf}
\begin{rem}\label{rem5.2}
In Examples \ref{exm5.1}, \ref
{exm5.2} and \ref{exm5.3}, the transition rates are
\textit{unbounded}, and the reward and costs are allowed to be \textit
{unbounded from above and from below}. In contrast, the transition rates
in \cite{havis,K85,K71,LP01,P94,S99,Y77} are assumed to be \textit
{bounded}, and the costs in \cite{G07a,OJ00,OJ03,OL99,K85,KNH}
are assumed to be \textit{nonnegative}. Moreover, Examples
\ref{exm5.2} and \ref{exm5.3} seem to be first computable
examples for the unconstrained optimal policies for discounted
continuous-time MDPs in Polish spaces.
\end{rem}

\section{Proofs of the main results}\label{proof}

In this section, we give proofs of Theorems \ref{t2},
\ref{Th-4.1}, \ref{Th-4.6}, \ref{Th-5.1}, \ref{Th-3.2},
\ref{Th-5.2} and of Lemmas \ref{Th-4.5} and \ref{Th-4.4}, which
are stated in Section \ref{section-3}.

To prove Theorems \ref{t2}, we need the following two lemmas.
\begin{lem}\label{Le-3.1}
Suppose that real-valued measurable
functions $\bar w\ge0$ on $S$ and $\bar q_t(D|x)$ on $
\R_+^0\times{\mathcal{B}}(S)\times S$ satisfy the following: for each
$t\geq0, D\in{\mathcal{B}}(S)$ and $x\in S$:
\begin{longlist}[(2)]
\item[(1)] $\bar q_t(\cdot|x)$ is a signed measure on
${\mathcal{B}}(S)$ such that $\bar q_t(S|x)\equiv0$, $\bar
q_t(D|x)\geq0$ for all $x\notin D$ and $\bar q_t(x):=\bar q_t(S-\{
x\}|x)<\infty$;\vspace*{1pt}
\item[(2)] $\int_S\bar w(y)\bar q_t(dy|x)\le\bar\rho\bar w(x)+\bar
b$, with constants $\bar\rho\ne0$ and $\bar b\ge0$.
\end{longlist}
Then nonnegative function
%
%
\begin{equation}\label{estar}
\bar h(s,x,t):=e^{\bar\rho(t-s)} \bar w(x)+\frac{\bar b}{\bar\rho
}\bigl(e^{\bar\rho(t-s)}-1\bigr)
\end{equation}
satisfies the following inequality:
\[
\int_s^t\int_{S-\{x\}}e^{-\int_s^z \bar q_v(x)\,dv}\bar q_z(dy|x)\bar
h(z,y,t)\,dz+e^{-\int_s^t \bar q_v(x)\,dv}\bar
w(x)\le\bar h(s,x,t)
\]
for all $x\in S$ and $0\leq s\le t<\infty$.
\end{lem}
\begin{pf}
Under conditions $(1)$ and $(2)$, a
straightforward calculation gives
\begin{eqnarray*}
&& \int_s^t\int_{S-\{x\}}e^{-\int_s^z \bar q_v(x)\,dv}\bar
q_z(dy|x)\bar
h(z,y,t)\,dz\\
&&\qquad \le\int_s^t e^{-\int_s^z \bar q_v(x)\,dv}\biggl[e^{\bar
\rho(t-z)}\biggl(\bar\rho\bar w(x)+\bar b\\
&&\qquad\hspace*{117.5pt}{}+\bar w(x)\bar
q_z(x)+\frac{\bar b}{\bar\rho} \bar q_z(x)\biggr)-\frac{\bar b}{\bar
\rho}\bar q_z(x)
\biggr]\,dz \\
&&\qquad =\bar h(s,x,t)-e^{-\int_s^t \bar q_v(x)\,dv}\bar w(x),
\end{eqnarray*}
which verifies this lemma.
\end{pf}
\begin{lem}\label{Le-3.2}
Suppose that Assumption \ref{assumA}(1) holds for
$\rho\ne0$. Then, for any $\pi\in\Pi$ and $x\in
S$,
\[
E_x^\pi\bigl[w(\xi_t)I_{\{t<T_{k+1}\}}\bigr]\le e^{\rho t}w(x)+\frac{b}{\rho
}(e^{\rho t}-1) \qquad\forall k\ge0 \mbox{ and } t\geq0,
\]
where $w$ and $b$ are from Assumption \ref{assumA}(1).
\end{lem}
\begin{pf}
Fix any $\pi\in\Pi, l\geq1$, and
$(x_0,\theta_1,x_1, \ldots,x_{l-1},\theta_l)\in(S\times
\R_+^0)^l$. Let $m_l(\cdot|h_l,t)$ be as in (\ref{3.1}). Then,
it follows from Assumption \ref{assumA}(1) that the following function
on $
\R_+^0\times{\mathcal{B}}(S)\times S$:
\[
\bar q_t(D|x):=\cases{
m_l(D|x_0,\theta_1,x_1,\ldots,\theta_l,x,t), &\quad if $x\notin D$,
\cr
-m_l(S|x_0,\theta_1,x_1,\ldots,\theta_l,x,t), &\quad if $D=\{x\}$,}
\]
satisfies conditions $(1)$ and $(2)$ for Lemma \ref{Le-3.1}.

Let $h(s,x,t):=e^{\rho(t-s)}w(x)+\frac{b}{\rho}(e^{\rho(t-s)}-1)$
for all $x\in S$ and $t\geq s\geq0$. Then, for each fixed $x\in
S$ and $0\leq s\leq t$, by Lemma \ref{Le-3.1} we have
%
%
\begin{eqnarray}\label{3.4}
&&
\int_s^t\int_{S-\{x\}}m_l(dy|h_{l-1},\theta
_l,x,z-T_l)h(z,y,t)\nonumber\\
&&\qquad\quad\hspace*{8.2pt}{}\times e^{-\int_s^z
m_l(S|h_{l-1},\theta_l,x,v-T_l)\,dv}\,dz
\nonumber\\
&&\quad{}
+w(x)e^{-\int_s^tm_l(S|h_{l-1},\theta_l, x,v-T_l)\,dv} \nonumber\\
&&\qquad=\int_{s-T_l}^{t-T_l}\int_{S-\{x\}}m_l(dy|h_{l-1},\theta
_l,x,\tilde
u)h(\tilde u,y,t-T_l)\\
&&\qquad\quad\hspace*{55.2pt}{}\times e^{-\int_{s-T_l}^{\tilde u}
m_l(S|h_{l-1},\theta_l,x,\tilde v)\,d\tilde v}\,d\tilde u \nonumber\\
&&\qquad\quad{}
+w(x)e^{-\int_{s-T_l}^{t-T_l} m_l(S|h_{l-1},\theta_l, x,\tilde
v)\,d\tilde v} \nonumber\\
&&\qquad \le h(s-T_l,x,t-T_l)=h(0,x,t-s).\nonumber
\end{eqnarray}
Moreover, by (\ref{2.4}) and (\ref{3.2}), we have
\begin{eqnarray*}
&&
E_x^\pi\bigl[w(\xi_t)I_{\{t<T_{k+1}\}}|\mathcal{F}_{T_k}\bigr]\\
&&\qquad=e^{-\int
_0^{t-T_k}m_k(S|h_k,v)\,dv}w(x_k)I_{\{T_k\le t\}}
+I_{\{T_k>t\}}\sum_{m=1}^k I_{\{T_{m-1}\le t<T_m\}}w(x_{m-1}).
\end{eqnarray*}
Now, using (\ref{3.4}) at $l=k,s=T_k=T_l,x=x_k=x_l$, gives
\begin{eqnarray*}
&&E^\pi_x\bigl[w(\xi_t)I_{\{t<T_{k+1}\}}|\mathcal{F}_{T_k}\bigr]\\
&&\qquad\le I_{\{T_k\le
t\}
}h(T_k,x_k,t)+I_{\{T_k>t\}}\sum_{m=1}^k
I_{\{T_{m-1}\le t<T_m\}}w(x_{m-1}),
\end{eqnarray*}
which implies that the following (\ref{3.5}) holds for $n=0$:
%
%
\begin{eqnarray} \label{3.5}
&&E^\pi_x\bigl[w(\xi_t)I_{\{t<T_{k+1}\}}|\mathcal{F}_{T_{k-n}}\bigr]\nonumber\\[-1pt]
&&\qquad\le I_{\{
T_{k-n}\le t\}}h(T_{k-n},x_{k-n},t)
+I_{\{T_{k-n}>t\}}\sum_{m=1}^{k-n} I_{\{T_{m-1}\le t<T_m\}}
w(x_{m-1}) \\[-1pt]
&&\eqntext{\forall k\geq n\geq0.}
\end{eqnarray}
Suppose that (\ref{3.5}) holds for some
$0\le n<k$. Then, by (\ref{3.2}) we have
\begin{eqnarray*}
&&E_x^\pi\bigl[w(\xi_t)I_{\{t<T_{k+1}\}}|\mathcal{F}_{T_{k-n-1}}\bigr]\\[-1pt]
&&\qquad \le E^\pi_x\Biggl[I_{\{T_{k-n}\le t\}}h(T_{k-n},x_{k-n},t)\\[-1pt]
&&\qquad\quad\phantom{ E^\pi_x\Biggl[}{} +I_{\{
T_{k-n}>t\}}\sum_{m=1}^{k-n}I_{\{T_{m-1}\le
t<T_m\}}w(x_{m-1})\Big|\mathcal{ F}_{T_{k-n-1}}\Biggr]\\[-1pt]
&&\qquad =E^\pi_x\bigl[I_{\{T_{k-n}\le t\}}h(T_{k-n},x_{k-n},t)\\[-1pt]
&&\qquad\quad\phantom{E^\pi_x\bigl[}{} +I_{\{
T_{k-n}>t\}} I_{\{T_{k-n-1}\le t<T_{k-n}\}}
w(x_{k-n-1})|\mathcal{F}_{T_{k-n-1}}\bigr] \\[-1pt]
&&\qquad\quad{} +I_{\{T_{k-n}>t\}}\sum_{m=1}^{k-n-1}I_{\{T_{m-1}\le
t<T_m\}}w(x_{m-1}) \\[-1pt]
&&\qquad =I_{\{T_{k-n-1}\le t\}}\biggl[\int_0^{t-T_{k-n-1}}\int_{S-\{
x_{k-n-1}\}}m_{k-n-1}(dy|h_{k-n-1},\tilde t\,)\\[-1pt]
&&\qquad\quad\hspace*{145.5pt}{}\times h(T_{k-n-1}+\tilde t,y,t)\\[-1pt]
&&\qquad\quad\hspace*{145.5pt}{}\times e^{-\int_0^{\tilde
t}m_{k-n-1}(S|h_{k-n-1},\tilde v)\,d\tilde v} \,d\tilde t\\[-1pt]
&&\qquad\quad\hspace*{81.5pt}{} +
e^{-\int_0^{t-T_{k-n-1}}m_{k-n-1}(S|h_{k-n-1},\tilde v)\,d\tilde
v}w(x_{k-n-1})\biggr] \\[-1pt]
&&\qquad\quad{} +I_{\{T_{k-n-1}>t\}}\sum_{m=1}^{k-n-1} I_{\{T_{m-1}\le
t<T_m\}}w(x_{m-1}),
\end{eqnarray*}
which together with $h(T_{k-n-1}+\tilde t,y,t)=h(\tilde
t,y,t-T_{k-n-1})$ and (\ref{3.4}) again, gives
\begin{eqnarray*}
&&E^\pi_x\bigl[w(\xi_t)I_{\{t<T_{k+1}\}}|\mathcal{F}_{T_{k-n-1}}\bigr]\\
&&\qquad\le I_{\{T_{k-n-1\}}\le
t\}}h(T_{k-n-1},x_{k-n-1},t)\\[-1pt]
&&\qquad\quad{}+I_{\{T_{k-n-1}>t\}}\sum_{m=1}^{k-n-1}
I_{\{T_{m-1}\le t<T_m\}}w(x_{k-1}).
\end{eqnarray*}
Hence, (\ref{3.5}) holds for all
$0\le n\le k$, and so this lemma follows from (\ref{3.5}) at
$n=k$.
\end{pf}
\begin{pf*}{Proof of Theorem \ref{t2}}
(a) We first prove the
following fact:
%
%
\begin{equation}\label{3.6}
P_x^\pi\bigl(\xi_tI_{\{T_k\le t<T_{k+1}\}}\notin S_l\mbox{:  for some }
k \geq0\bigr)\to0 \qquad\mbox{as } l\to\infty.
\end{equation}
To prove\vspace*{1pt} (\ref{3.6}), let $\Gamma_l:=\{e\dvtx\xi_t(e)I_{\{T_k\le
t<T_{k+1}\}}(e)\notin S_l$ for some $k\geq0\}$ for any
$l\geq1$.

Suppose that, for some $\varepsilon>0$ and any $\tilde L\ge1$,
there exists $l>\tilde L$ such that
%
%
\begin{equation}\label{3.7}
P^\pi_x(\Gamma_l)= P^\pi_x\bigl(\bigl\{e\dvtx\xi_t(e)I_{\{T_k\le
t<T_{k+1}\}}(e)\notin S_l \mbox{ for some } k\geq0
\bigr\}\bigr)>\varepsilon.
\end{equation}
Then, by Assumption \ref{assumA}(2), we can take the corresponding $l$ such
that (\ref{3.7}) holds and also
the following inequality:
%
%
\begin{equation}\label{3.8}
w(y)>\biggl[e^{\tilde\rho t}w(x)+\frac{b}{\tilde\rho}(e^{\tilde\rho
t}-1)\biggr]\Big/\varepsilon\qquad\forall y\notin S_l,
\end{equation}
is satisfied, where $\tilde\rho:=|\rho|+1$.

For the taken $l\geq1$ in (\ref{3.8}), let us define new
transition rates $\tilde q(D|x,a)$ as follows:
\[
\tilde q(D|x,a):=\cases{
q(D|x,a), &\quad if $x\in S_l$, \cr
0, &\quad if $x\notin S_l$,}
\qquad\mbox{for }
(x,a)\in K.
\]
The
quantities such as probabilities corresponding to $\tilde q(D|x,a)$
are equipped with the tilde.

We next to prove that
%
%
\begin{eqnarray} \label{3.9}
&&
P^\pi_x\bigl(\xi_tI_{\{T_k\le t<T_{k+1}\}}\in S_l \mbox{ for all }
k\ge0\bigr)\nonumber\\[-8pt]\\[-8pt]
&&\qquad=\tilde P^\pi_x\bigl(\xi_tI_{\{T_k\le t<T_{k+1}\}}\in S_l
\mbox{ for all } k \ge0\bigr).\nonumber
\end{eqnarray}
Indeed, it is obvious that
\[
P^\pi_x(X_0\in S_l)=\tilde P^\pi_x(X_0\in S_l)=I_{{S_l}}(x).
\]
Let $X_k^t:=X_kI_{\{T_k\le t<T_{k+1}\}}$. Then, by (\ref{2.4}) we
have $\{\xi_tI_{\{T_k\le t<T_{k+1}\}}\in S_l\}=\{X_k^t\in S_l\}$.
We now suppose that for some $n\geq0$,
%
%
\begin{eqnarray}\label{3.10}
&&
P^\pi_x(\{X_k^t\in S_l, 0\leq k\leq n\}\cap\Gamma)\nonumber\\[-8pt]\\[-8pt]
&&\qquad=\tilde
P^\pi_x( \{X_k^t\in S_l, 0\leq k\leq n\}\cap\Gamma) \qquad
\forall\Gamma\in{\mathcal{B}}(\hat H_n),\nonumber
\end{eqnarray}
where $P_x^\pi$ and $\tilde P_x^\pi$ are regarded as the marginal on
$\hat H_{n+1}$.

Using the notation in (\ref{3.1}) and (\ref{3.2}), for any $D\in
{\mathcal{B}}(S),
0<t_1<t_2<\infty$, we have
\begin{eqnarray*}
&&P^\pi_x\bigl(\{X_k^t\in S_l, 0\leq k\leq n\mbox{,  and } X_{n+1}^t\in
S_l\}\cap
\{\Gamma\times(t_1,t_2)\times D\}\bigr) \\
&&\qquad=\int_{{t_1}}^{{t_2}}\int_\Gamma P^\pi_x(dh_n)I_{\{X_k^t\in S_l,
0\leq k\leq n\}}I_{\{X_{n+1}^t\in S_l\cap D\}}\\
&&\qquad\quad\hspace*{28pt}{}\times m_n(S_l\cap D|h_n,\tilde t\,)e^{-\int_0^{{\tilde t}}
m_n(S|h_n,v)\,dv}\,d\tilde t \\
&&\qquad=\int_{{t_1}}^{{t_2}}\int_\Gamma\tilde
P^\pi_x(dh_n)I_{\{X_k^t\in S_l, 0\leq k\leq n\}}I_{\{X_{n+1}^t\in
S_l\cap D\}} \\
&&\qquad\quad\hspace*{28pt}{}\times \tilde m_n(S_l\cap D|h_n,\tilde t\,)e^{-\int_0^{{\tilde
t}} \tilde m_n(S|h_n,v)\,dv}\,d\tilde t
\\
&&\qquad=\tilde P^\pi_x\bigl(\{X_k^t\in S_l, 0\leq k\leq n\mbox{,  and}
X_{n+1}^t\in S_l\}\cap
\{\Gamma\times(t_1,t_2)\times D\}\bigr),
\end{eqnarray*}
which together with the arbitrariness of $D\in{\mathcal{B}}(S)$ and
$0\leq t_1<t_2$ implies (\ref{3.10}) for $n+1$, and thus
(\ref{3.9}) follows from the induction.

Thus, from (\ref{3.7}) and (\ref{3.9}), we have
%
%
\begin{equation}\label{3.11}
\tilde P^\pi_x(\Gamma_l)=\tilde P^\pi_x\bigl(\xi_tI_{\{T_k\le t<T_{k+1}\}
}\notin S_l \mbox{ for some } k\geq0
\bigr)>\varepsilon.
\end{equation}
Moreover, since $\|\tilde q\|:=\sup_{x\in S, a\in A(x)}|\tilde
q(\{x\}|x,a)|=\sup_{x\in S_l, a\in A(x)}|q(\{x\}|x,\break a)|<\infty$, we
now show by induction that
%
%
\begin{equation}\label{3.12}
\tilde E_x^\pi[e^{-T_k}]\leq\bigl[1-e^{-\|\tilde q\|}(1-e^{-1})\bigr]^k
\qquad\forall k\geq1.
\end{equation}
In fact, by (\ref{3.1}) we have $|\tilde m_k(S|h_k)\leq
\|\tilde q\|$ for all $k\geq1$, and it follows from (\ref{3.2})
that
%
%
\begin{eqnarray} \label{3.13}\qquad
\tilde E_x^\pi[e^{-T_1}]&=&\int_0^1\tilde
m_0(S|x)e^{-\tilde m_0(S|x)t}e^{-t}\,dt+\int_1^\infty\tilde
m_0(S|x)e^{-\tilde
m_0(S|x)t}e^{-t}\,dt \nonumber\\[-8pt]\\[-8pt]
&\leq& 1-e^{-\|\tilde q\|}\int_0^1e^{-t}\,dt=\bigl[1-e^{-\|\tilde q\|
}(1-e^{-1})\bigr].\nonumber
\end{eqnarray}
Suppose that (\ref{3.12}) holds for some $k\geq1$. Then, as the
arguments of (\ref{3.13}), from (\ref{3.1}) and (\ref{3.2}) we
also have $\tilde E_x^\pi[e^{-T_{k+1}}]\leq\tilde
E_x^\pi[e^{-T_k} [1-e^{-\|\tilde
q\|}(1-e^{-1})]]\leq[1-e^{-\|\tilde q\|}(1-e^{-1})]^{k+1}$, and so
(\ref{3.12}) follows. Hence, by (\ref{3.12}) and the Chebychev
inequality we have
\begin{eqnarray*}
\tilde P_x^\pi(T_\infty\leq t)&\leq&\tilde P_x^\pi(T_k \leq t)=\tilde
P_x^\pi(e^{-T_k} \geq e^{-t})\leq e^t\tilde
E_x^\pi[e^{-T_k}]\\
&\leq& e^t\bigl[1-e^{-\|\tilde q\|}(1-e^{-1})\bigr]^k
\end{eqnarray*}
for all $k\geq1$, and so $\tilde P_x^\pi(T_\infty\geq t)=1$.
Since $t>0$ can be arbitrary, we have
$\tilde P_x^\pi(T_\infty=\infty)=1$, and therefore, $\sum
_{k=0}^\infty
\tilde P^\pi_x(T_k\le t<T_{k+1})=1$. Since Assumption \ref{assumA}(1)
still holds
when $\rho$ and $q(D|x,a)$ are replaced with $\bar\rho$ and $\tilde
q(D|x,a)$, respectively, by Lemma \ref{Le-3.2} we have
%
%
\begin{equation} \label{3.14}\quad
\tilde E_x^\pi[w(\xi_t)] =\lim_{k\to\infty} \tilde
E_x^\pi\bigl[w(\xi_t)I_{\{t< T_{k+1}\}}\bigr]\le e^{\tilde\rho
t}w(x)+\frac{b}{\tilde\rho}(e^{\tilde\rho t}-1).
\end{equation}
On the other hand, using (\ref{3.8}) and (\ref{3.11}), we see
\begin{eqnarray*}
\tilde E_x^\pi[w(\xi_t)]&=&\tilde E_x^\pi[w(\xi_t)|\Gamma_l]\tilde
P_x^\pi(\Gamma_l)+ \tilde E_x^\pi[w(\xi_t)|\Gamma_l^c]\tilde
P_x^\pi(\Gamma_l^c)\\
&>&e^{\tilde\rho
t}w(x)+\frac{b}{\tilde\rho}(e^{\tilde\rho t}-1),
\end{eqnarray*}
which
contradicts to (\ref{3.14}), and thus (\ref{3.6}) is proved.

Since $\Gamma_{l+1}\subseteq\Gamma_l$ for all $l\geq1$, by (\ref
{3.6}) we conclude
that $P^\pi_x(\bigcap_{l\ge0}\Gamma_l)=0$, and so
%
%
\begin{equation}\label{3.15}\qquad
P^\pi_x\bigl(\bigl\{ \mbox{for each } l\geq1\mbox{,  there exists } k
\mbox{ such that } \xi_tI_{\{T_k\le
t<T_{k+1}\}}\notin S_l\bigr\}\bigr)=0.
\end{equation}
Since $\{\inf\{s\dvtx\xi_s \notin S_l\}\le t\}
\subseteq\{\xi_tI_{\{T_k\le t<T_{k+1}\}} \notin S_l$, for
some $k\geq1\}$, by (\ref{3.15}) we conclude $P^\pi_x(\inf\{s\dvtx
\xi_s \notin S_l\}\le t, l=1,\ldots)=0$, and thus $P^\pi_x(
\inf\{s\dvtx\xi_s \notin S_l\}>t$, for some $l\geq1)=1$,
or, equivalently, $P^\pi_x(\xi_s\in S_l$ for all $s\in[0,t]$,
for some $l\geq1)=1$. For any $k\geq1$, let
$B_k:=\{\xi_s\in S_l$ for all $s\in[0,k]$, for
some $l\geq1\}$. Then, $B_{k+1}\subseteq B_k$ and $P_x^\pi(B_k)=1$ for
all $k\geq1$, and thus $P^\pi_x(\bigcap_{k=1}^\infty B_k)=1$, which
together with (\ref{2.4}) implies $P^\pi_x(T_\infty=\infty)=1$. To
further prove \mbox{$P^\pi_x(\xi_t\in S)=1$}, using the facts $\sum_{k\ge
0} P^\pi_x(T_k\le t<T_{k+1})=P_x^\pi(T_\infty=\infty)=1$ and $
P^\pi_x(\xi_t\in S|T_k\le t<T_{k+1})=1$ for all $k\geq1$, we have
that $P^\pi_x(\xi_t\in S)=\sum_{k\ge0} P^\pi_x(\xi_t\in S|T_k\le
t<T_{k+1})P^\pi_x(T_k\le
t<T_{k+1})=1$, and thus (a) follows.

(b) First, consider the case of $\rho\ne0$. Since
$\sum_{k=0}^\infty P^\pi_x(T_k\le t<T_{k+1})=1$ for all $t\geq0$,
\[
E^\pi_x[w(\xi_t)]=E^\pi_i\Biggl[w(\xi_t)\sum_{k=0}^\infty I\{
T_k\le
t<T_{k+1}\}\Biggr]=\lim_{k\to\infty} E^\pi_i[
w(\xi_t)I\{t<T_{k+1}\}],
\]
which together with Lemma \ref{Le-3.2} implies the first part of
(b). Moreover, the results for the case of $\rho=0$ can be
obtained by letting $\rho\downarrow0$.

(c) Define an integer-valued random measure $\tilde\mu^*$ on
${\mathcal{B}}(\R_+^0)\times{\mathcal{B}}(S)$
%
%
\begin{equation}\label{4.4}
\tilde\mu^*(dt,dx):=\sum_{k\ge
1}I_{\{T_k<\infty\}}\delta_{(T_k,X_{k-1})}(dt,dx),
\end{equation}
which counts the exits from $dx$. Then, as Lemma 4.28 in
\cite{KR95}, the random measure
\[
\tilde\nu^\pi(e,dt,dx):=-\biggl[\int_A
\pi(da|e,t)q(dx|\xi_{t-}(e),a)I_{dx}(\xi_{t-}(e))\biggr]\,dt
\]
is a dual
predictable projection of the measure $\tilde\mu^*$ with respect
to $\mathcal P$ and $P^\pi_\gamma$ (for any fixed policy $\pi\in\Pi$
and initial distribution $\gamma$). Hence, by (4.5) in \cite{KR95}
we have
\begin{eqnarray*}
E^\pi_x[\tilde\mu^*((0,t],D)]&=&E^\pi_x[\tilde\nu^\pi((0,t],D)]\\
&\le&
E^\pi
_x\biggl[\int_0^t\int_A\pi(da|e,s)\sup_{x\in D}q^*(x)\,ds\biggr]
<\infty
\qquad\forall
t\geq0,
\end{eqnarray*}
which together with
$|\mu^*((0,t],D)-\tilde\mu^*((0,t],D)|\le1$ and (4.5) in
\cite{KR95} again, implies
\[
E^\pi_x[\mu^*((0,t],D)]=E^\pi_x[\nu^\pi((0,t],D)]<\infty.
\]
Thus, using the obvious representation
$I_{\{\xi_t\in D\}}=I_D(x)+\mu^*((0,t],D)-\tilde\mu^*((0,\break t],D)$, by
taking the expectation $E^\pi_x$ of the representation we see that (c)
is true.
\end{pf*}
\begin{pf*}{Proof of Theorem \ref{Th-4.1}} (a) For the given $D$, by
Theorem \ref{t2}(c) and (\ref{4.1}) we have
\begin{eqnarray*}
\hat\eta^\pi(D)&=&\gamma(D)+\alpha\int_0^\infty e^{-\alpha
t}E^\pi_\gamma\biggl[\int_0^t \int_A\pi(da|e,s)q(D|\xi
_{s-}(e),a)\,ds
\biggr]\,dt \\
&=& \gamma(D)+\alpha\int_S\int_{A} q(D|x,a) \\
&&\hspace*{68.5pt}{}\times\int_0^\infty
e^{-\alpha
t} \int_0^t
E^\pi_\gamma\bigl[\pi(da|e,s)I_{\{\xi_{s-}(e)\}}(dx)\,ds\bigr]\,dt \\
&=& \gamma(D)+\frac{1}{\alpha} \int_S\int_{A} q(D|x,a)\eta^\pi(dx,da),
\end{eqnarray*}
and so (a) follows.

(b) Recall that $\eta(dx,da)=\hat\eta(dx)\phi^\eta(da|x)$. Then,
to prove (b), it suffices to show
%
%
\begin{equation}\label{4.7}
\int_{S}\int_{A}
u(x,a)\eta(dx,da)=\int_{S}\int_Au(x,a)\eta^{\phi^\eta}(dx,da)
\end{equation}
for each nonnegative bounded measurable function $u$ on $K$. In fact,
for any such a function $u$, by Lemma 5.3 in \cite{G07b} and (\ref
{2.10}) we have
%
%
\begin{eqnarray}\label{4.8}
\alpha V_\alpha(x,\phi^\eta,u)&=&\int_{A(x)}
u(x,a)\phi^\eta(da|x)\nonumber\\[-8pt]\\[-8pt]
&&{}+
\int_S V_\alpha(y,\phi^\eta,u)q(dy|x,\phi^\eta) \qquad\forall
x\in S.\nonumber
\end{eqnarray}
On the other hand, let $\|u\|_1:=\sup_{(x,a)\in K}|u(x,a)|
<\infty$, and $|q(dx|x,\phi^\eta)|$ the total variation of
$q(dy|x,\phi^\eta)$. Then, by $(T_2)$--$(T_3)$ and the
condition in (b) we have
\[
\int_S\int_S|V_\alpha(y,\phi^\eta,u)||q(dy|x,\phi^\eta)|\hat
\eta(dx)
\leq\frac{2\|u\|_1}{\alpha} \int_S |q(\{x\}|x,\phi^\eta)|\hat
\eta
(dx) <\infty,
\]
which together with the Jordan decomposition of
$q(\cdot|x,\phi^\eta)$ and Theorem 2.6.4 in~\cite{A00}, implies
\[
\int_{S}
\int_S[\hat\eta(dy)q(dx|y,\phi^\eta)]V_\alpha(x,\phi^\eta
,u)=\int_S\biggl[\int
_SV_\alpha(y,\phi^\eta,u)q(dy|x,\phi^\eta)\biggr]\hat\eta(dx).
\]
Hence, by Assumption \ref{assumA}(3) we have
%
%
\begin{eqnarray}\label{q}
&&\lim_{k\to\infty}\int_{S_k}
\int_S[\hat\eta(dy)q(dx|y,\phi^\eta)]V_\alpha(x,\phi^\eta
,u)\nonumber\\[-8pt]\\[-8pt]
&&\qquad=\lim_{k\to
\infty}\int_{S_k}
\biggl[\int_SV_\alpha(y,\phi^\eta,u)q(dy|x,\phi^\eta)\biggr]\hat\eta(dx).\nonumber
\end{eqnarray}
Thus, for any fixed $k\geq1$, since $\sup_{x\in
S_k}q^*(x)<\infty$, by (\ref{4.8}) and (a) we have
\begin{eqnarray*}
&&
\int_{S_k}\int_{A} u(x,a)\eta(dx,da)\\
&&\qquad=\int_{S_k}\int_{A(x)}
u(x,a)[\hat
\eta(dx)\phi^\eta(da|x)] \\
&&\qquad=\int_{S_k}\biggl[\alpha V_\alpha(x,\phi^\eta,u)-\int_{S}V_\alpha
(y,\phi
^\eta,u)q(dy|x,\phi^\eta)\biggr]\hat\eta(dx)\\
&&\qquad=\alpha\int_{S_k}V_\alpha(x,\phi^\eta,u) \gamma(dx)+\int_{S_k}
V_\alpha(y,\phi^\eta,u)\biggl[\int_{S}\hat\eta(dx)q(dy|x,\phi^\eta)\biggr]
\\
&&\qquad\quad{} -\int_{S_k}\biggl[\int_{S}V_\alpha(y,\phi^\eta,u)q(dy|x,\phi
^\eta
)\biggr]\hat\eta(dx)\\
&&\qquad=\int_{S_k}\int_A u(x,a)\eta^{\phi^\eta}(dx,da)+\int_{S_k}
\biggl[\int_{S}\hat\eta(dy)q(dx|y,\phi^\eta)\biggr]V_\alpha(x,\phi^\eta,u)
\\
&&\qquad\quad{}
-\int_{S_k}\biggl[\int_{S}V_\alpha(y,\phi^\eta,u)q(dy|x,\phi^\eta
)\biggr]\hat\eta(dx),
\end{eqnarray*}
which together with (\ref{q}) gives (\ref{4.7}).

(c) Since $\phi\in\Pi_s$, by (a) and (\ref{4.2}) we have
\begin{eqnarray}
\alpha\hat\eta^\phi(D)&=&\alpha\gamma(D)+\int_Sq(D|x,\phi)\hat
\eta
^\phi(dx)\nonumber\\
&=&\alpha\gamma(D)+\int_S\int_Aq(D|x,a)[\hat\eta^\phi(dx)\phi(da|x)]
\nonumber\\
&&\eqntext{\forall D\in{\mathcal{B}}(S) \mbox{ with } \displaystyle \sup_{x\in D}
q^*(x)<\infty.}
\end{eqnarray}
Moreover, under Assumptions \ref{assumA}, \ref{assumB}(2) and \ref
{assumB}(3), by Theorem \ref{t3}
we have
%
%
\begin{equation} \label{M^*}\qquad
\int_S|q(\{x\}|x,\phi)|\hat\eta^\phi(dx)\leq L\biggl[\alpha\int_S
w(x)\gamma(dx)+b\biggr]\Big/[\alpha(\alpha-\rho)]<\infty.
\end{equation}
Thus, by (b) we see that
$\hat\eta^\phi(dx)\phi(da|x)=\eta^\phi(dx,da)$, and so (c)
follows.
\end{pf*}
\begin{pf*}{Proof of Lemma \ref{Th-4.5}} (a) Since the first part of (a)
follows from (\ref{4.13}), we need to verify the second part of
(a). In fact, for each $\mu\in{\mathcal{P}}(S\times A)$, by
(\ref{4.14}) we have $\int_S \bar w(x)\hat T_{\bar
w}'(\mu)(dx)=\frac{1}{\int_S{1}/({\bar
w(x)})\hat\mu(dx)}<\infty$, and so the second part of (a) follows.

(b) By (\ref{4.13}) and (\ref{4.14}) and a straightforward calculation,
we see that (b) is true.

(c) and (d). We prove (c) and (d) together. Suppose that $\eta_k
\stackrel{\bar w}{\longrightarrow}\eta_0$. Take any bounded
continuous function $u$ on $S\times A$. Then, since $\bar w$ is
continuous, by $\eta_k \stackrel{\bar w}{\longrightarrow}\eta_0$
we have
\begin{eqnarray*}
&&\lim_{k\to\infty}\int_{S\times A}v(x,a)\bar
w(x)\eta_k(dx,da)\\
&&\qquad=\int_{S\times A}v(x,a)\bar w(x)\eta_0(dx,da)
\qquad\mbox{for } v:=u, 1,
\end{eqnarray*}
which together with (\ref{4.13}), imply
%
%
\begin{equation}\label{4.15}\quad
\lim_{k\to\infty}\int_{S\times A}u(x,a)T_{\bar w}(\eta_k)(dx,da) =
\int_{S\times A}u(x,a)T_{\bar w}(\eta_0)(dx,da),
\end{equation}
and thus, $T_{\bar w}(\eta_k)\stackrel{1}{\longrightarrow}T_{\bar
w}(\eta_0)$.

On the other hand, suppose that $\mu_k
\stackrel{1}{\longrightarrow}\mu_0$, and pick up any continuous
function $u(x,a)$ on $S\times A$ such that $|u(x,a)|\leq L_u\bar
w(x)$ for all $(x,a)\in K$, with some nonnegative constant $L_u$
depending on $u$. Then, the functions $\frac{u(x,a)}{\bar w(x)}$
and $\frac{1}{\bar w}$ are bounded continuous on $S\times A$.
Hence, a straightforward calculation gives
%
%
\begin{equation} \label{4.16}\quad
\lim_{k\to\infty}\int_{S\times A}u(x,a)T_{\bar w}'(\mu_k)(dx,da)
=\int_{S\times A}u(x,a)T_{\bar w}'(\mu_0)(dx,da).
\end{equation}
By (\ref{4.15}) and (\ref{4.16}) and (b), we see that (c) and (d) are
both true.
\end{pf*}
\begin{pf*}{Proof of Lemma \ref{Th-4.4}} (a) For any $\eta^{\pi_1},
\eta^{\pi_2}\in{\mathcal{M}}_o$ and $0\leq\beta\leq1$, let
$\eta:=\beta\eta^{\pi_1}+(1-\beta)\eta^{\pi_2}$. Then, by
Theorem \ref{Th-4.1}(a) and a straightforward calculation we have
%
%
\begin{eqnarray}\label{4.10}
\alpha\hat\eta(D)&=&\alpha\gamma(D)+\int_{S\times
A}q(D|x,a)\eta(dx,da) \nonumber\\[-8pt]\\[-8pt]
&&\eqntext{\displaystyle \forall D\in{\mathcal{B}}(S) \mbox{ with
} \sup_{x\in D} q^*(x)<\infty,}
\end{eqnarray}
and also
$\int_Sw(x)\hat\eta(dx)=\int_Sw(x)[\beta\hat\eta^{\pi
_1}(dx)+(1-\beta
)\hat\eta^{\pi_2}(dx)]<\infty$.
Thus, by Theorem \ref{Th-4.1}(b) and (\ref{4.10}), there exists a
randomized stationary policy $\phi^\eta\in\Pi_s$ such that
$\eta=\eta^{\phi^\eta}$. Hence, ${\mathcal{M}}_o$ is convex, and thus
so is ${\mathcal{M}}_o^c$.

(b) Take any sequence $\{\eta_m\}$ in ${\mathcal{M}}_o$ such
that $\eta_m \stackrel{w}{\longrightarrow} \eta_0$ (and thus
$\eta_m \stackrel{1}{\longrightarrow} \eta_0$). Then, under
Assumptions \ref{assumA}, \ref{assumB}(2) and \ref{assumB}(3), by
Theorem \ref{t2}(b) we have
%
%
\begin{eqnarray} \label{G-1}
\int_Sw(x)\hat\eta_m(dx)&=&\int_Sw(x)\eta_m(dx,da)\leq
\frac{\alpha\int_Sw(x)\gamma(dx)+b}{\alpha(\alpha-\rho
)}\nonumber\\[-8pt]\\[-8pt]
&=&M^*_1<\infty
\qquad\forall m\geq1.\nonumber
\end{eqnarray}
Thus, by Lemma 11.4.7 in \cite{OL99} we have
\begin{eqnarray*}
\int_S|q(\{x\}|x,\phi^{\eta_0})|\hat\eta_0(dx)&\leq& L\int
_Sw(x)\hat\eta
_0(dx)\leq L\liminf_{m\to\infty}\int_Sw(x)\hat\eta_m(dx)\\
&\leq&
LM^*_1<\infty.
\end{eqnarray*}
Thus, to prove $\eta_0\in{\mathcal{M}}_o$, by Theorem
\ref{Th-4.1}(b) it suffices to show
\begin{eqnarray}
\alpha\hat\eta_0(D)&=&\alpha\gamma(D)+\int_K q(D|x,a)\eta_0(dx,da)
\nonumber\\
&&\eqntext{\forall D\in{\mathcal{B}}(S) \mbox{ with } \displaystyle \sup_{x\in D}
q^*(x)<\infty,}
\end{eqnarray}
which can follow (by Proposition 7.18 in \cite{b7}) from
%
%
\begin{eqnarray}\label{4.11}\qquad
\alpha\int_S g(y)\hat\eta_0(dy)
=\alpha
\int_Sg(y)\gamma(dy)
+\int_S \int_K g(y)q(dy|x,a)\eta_0(dx,da)
\nonumber\\[-8pt]\\[-8pt]
&&\eqntext{\forall g\in C_b(S).}
\end{eqnarray}
Thus, the rest verifies (\ref{4.11}). For any $g\in C_b(S)$, by
$\eta_m\in{\mathcal{M}}_o$ and Theorem \ref{Th-4.1}(a) we have
%
%
\begin{eqnarray}\label{G-2}\hspace*{32pt}
\alpha\int_{S_k}g(y)\hat\eta_m(dy)=\alpha
\int_{S_k}g(y)\gamma(dy)+\int_{S_k}\int_K
g(y)q(dy|x,a)\eta_m(dx,da) \nonumber\\[-8pt]\\[-8pt]
&&\eqntext{\forall k, m\geq1.}
\end{eqnarray}
Since $q^*(x)\leq Lw(x)$ for all $x\in S$, using Assumption \ref{assumA}(3)
and the dominated convergence theorem, by (\ref{G-1}) and
(\ref{G-2}) with letting $k\to\infty$ we have
%
%
\begin{eqnarray}\label{4.12}\hspace*{32pt}
\alpha\int_S g(y)\hat\eta_m(dy)=\alpha
\int_Sg(y)\gamma(dy)+\int_S \int_K
g(y)q(dy|x,a)\eta_m(dx,da)\nonumber\\[-8pt]\\[-8pt]
&&\eqntext{\forall m\geq1.}
\end{eqnarray}
On the other hand, since $|\int_S g(y)q(dy|x,a)|\leq2\|g\|_1
q^*(x)\leq2L\|g\|_1w(x)$ [for all $a\in A(x)$], by $\eta_m
\stackrel{w}{\longrightarrow} \eta_0$ and Assumption \ref{assumC}(1),
we have
\begin{eqnarray*}
\lim_{m\to\infty} \int_S g(y)\hat\eta_m(dy)&=&\lim_{m\to\infty}
\int_S g(y)\eta_m(dy,da)=\int_S
g(y)\eta_0(dy,da)\\
&=&\int_S g(y)\hat\eta_0(dy)
\end{eqnarray*}
and
\[
\lim_{m\to\infty}\biggl[\int_S\int_K
g(y)q(dy|x,a)\eta_m(dx,da)\biggr]=\int_S \int_K
g(y)q(dy|x,a)\eta_0(dx,da),
\]
which together with (
\ref{4.12}) give (\ref{4.11}), and so (b) follows.
\end{pf*}
\begin{pf*}{Proof of Theorem \ref{Th-4.6}} (a) Since ${\mathcal
{P}}(S\times
A)$ is metrizable, it follows from Lemma \ref{Th-4.5} (with $\bar
w:=w$) that ${\mathcal{P}}_{w}(S\times A)$ is also metrizable, and so
are ${\mathcal{M}}_o$ and ${\mathcal{M}}_o^c$. Since ${\mathcal
{M}}_o$ is
closed (by Lemma \ref{Th-4.4}) and ${\mathcal{M}}_o^c$ is a closed
subset of ${\mathcal{M}}_o$ under the additional Assumption \ref
{assumC}(1), it
suffices to show that ${\mathcal{M}}_o$ is sequentially relatively
compact. Indeed, for each $\eta\in{\mathcal{M}}_o$, since
$1\leq\int_S w'(x)\hat\eta(dx)<\infty$ [using Assumption \ref{assumC}(2)],
$T_{w'}(\eta)$ is well defined. Moreover, by (\ref{4.13}) and
Theorem \ref{t3}, we have
\begin{eqnarray*}
\int_{S\times
A}\frac{w(x)}{w'(x)}T_{w'}(\eta)(dx,da)&=&\frac{\int_{S\times
A}w(x)\eta(dx,da)}{\int_{S\times
A}w'(x)\eta(dx,da)}\\
&\leq&\int_{S\times A} w(x)\eta(dx,da)\leq\alpha
M_1^* \qquad\forall\eta\in{\mathcal{M}}_o,
\end{eqnarray*}
where $M_1^*$ is as in Theorem \ref{t3}(b). Thus, by Assumption
\ref{assumC}(2) and Prohorov' theorem (see Theorem 12.2.15 in
\cite{OL99}) we see that $\{T_{w'}(\eta), \eta\in{\mathcal{M}}_o\}$
is sequentially relatively compact, and so is ${\mathcal{M}}_o$ (by
Lemma \ref{Th-4.5} with $\bar w:=w'$).

(b) Under Assumptions \ref{assumA} and \ref{assumB}, by Theorem \ref
{t3}(b) we have
$|V_r(\pi)|\leq MM^*_1$ and $|V_n(\pi)|\leq MM^*_1$ for $1\leq
n\leq N$. Moreover, by Theorem \ref{Th-4.1} and (\ref{P})
[equivalently, (\ref{4.2b})] we can find a sequence
$\{\eta^{\pi_k}\}$ $(\pi_k\in\Pi_s, k=1,\ldots)$ such that
%
%
\begin{eqnarray}\label{4.17}
&\displaystyle V_r(U)=\lim_{k\to\infty}\frac{1}{\alpha}\int_K
r(x,a)\eta^{\pi_k}(dx,da),& \nonumber\\[-8pt]\\[-8pt]
&\displaystyle \int_K
c_n(x,a)\eta^{\pi_k}(dx,da)\leq\alpha d_n,\qquad n=1,\ldots, N.&\nonumber
\end{eqnarray}
Then, by (a) there exists a subsequence $\{\eta^{\pi_{k_m}}\}$ and
$\eta
_0\in{\mathcal{M}}_o$ such that
$\eta^{\pi_{k_m}} \stackrel{w}{\longrightarrow} \eta_0$ as
$m\to\infty$, which together with (\ref{4.17}) implies
\[
V_r(U)=\frac{1}{\alpha}\int_K r(x,a)\eta_0(dx,da)
\]
and
\[
\int
_K c_n(x,a)\eta_0(dx,da)\leq\alpha d_n,\qquad n=1,\ldots, N,
\]
and so $\phi^{\eta_0}$ is constrained optimal.
\end{pf*}
\begin{pf*}{Proof of Theorem \ref{Th-5.1}} Obviously, parts (a), (b) are
directive consequence of (\ref{LP}) and Theorem \ref{Th-4.1}.
Moreover, (c) follows from (b) and Theorem~\ref{Th-4.6}(b).
\end{pf*}
\begin{pf*}{Proof of Theorem \ref{Th-3.2}} (a) Under Assumptions \ref{assumA},
\ref{assumB}(2), \ref{assumB}(3) and \ref{assumC}(3), by Theorems \ref
{t2} and \ref{Th-4.1} we have
\begin{eqnarray*}
{\mathcal{M}}_o&=&\biggl\{\eta^\pi\Big| \int_Sw(x)\hat\eta^\pi(dx)\leq
\alpha
M^*_1, \pi\in\Pi\biggr\}\\
&=&\biggl\{\eta^\pi\Big|
\int_Sw(x)\hat\eta^\pi(dx)\leq\alpha M^*_1, \pi\in\Pi_s\biggr\}.
\end{eqnarray*}

We now prove that $\eta^f$ is an extreme point in ${\mathcal{M}}_o$
for each $f\in F$. In fact, for any fixed $f\in F$, suppose that
$\eta^f$ is not any extreme in ${\mathcal{M}}_o$. Then, there exist
$\beta\in(0,1)$ and $\pi_1,\pi_2\in\Pi_s$ such that
%
%
\begin{equation}\label{4.19}
\eta^f=\beta\eta^{\pi_1}+(1-\beta)\eta^{\pi_2}\quad \mbox{and}\quad
\eta^{\pi_1}\not=\eta^{\pi_2},
\end{equation}
which implies that $\hat\eta^{\pi_k}\ll\hat\eta^f$ ($k=1,2$).
Thus, it follows from (\ref{4.19}) and Theorem~\ref{Th-4.1} that
%
%
\begin{eqnarray}\label{4.20}
&\displaystyle
f(da|x)=\beta\frac{d\hat\eta^{\pi_1}}{d\hat\eta^{f}}(x)\pi
_1(da|x)+(1-\beta)\frac{d\hat\eta^{\pi_2}}{d\hat\eta^{f}}(x)\pi_2(da|x)\quad
\mbox{and} &\nonumber\\[-8pt]\\[-8pt]
&\displaystyle  \beta\frac{d\hat\eta^{\pi_1}}{d\hat\eta^{f}}(x)+(1-\beta
)\frac
{d\hat\eta^{\pi_2}}{d\hat\eta^{f}}(x)=1 \qquad
\forall x\in\hat S&\nonumber
\end{eqnarray}
for some\vspace*{-1pt} $\hat S\in{\mathcal{B}}(S)$ with $\hat\eta^f(\hat S)=1$,
where $\frac{d\hat\eta^{\pi_k}}{d\hat\eta^{f}}$ denote the
(nonnegative) Radon--Nikodym derivative. Moreover, by
$\eta^{\pi_1}\not=\eta^{\pi_2}$ we see that $\hat\eta^f(\{x\in
\hat S| \pi_1(\Gamma|x)\not=\pi_2(\Gamma|x)$ for some
$\Gamma\in{\mathcal{B}}(A)\})>0$. (Otherwise, $\eta^{\pi_1}$ and
$\eta^{\pi_2}$ coincide.) Thus, for each $x\in\{x\in\hat S|
\pi_1(\Gamma|x)\not=\pi_2(\Gamma|x)$ for some
$\Gamma\in{\mathcal{B}}(A)\}$, there exists a corresponding
$\Gamma_x\in{\mathcal{B}}(A)$ (depending on $x$) such that $0<
\pi_1(\Gamma_x|x)<\pi_2(\Gamma_x|x)<1$. Therefore, by (\ref{4.20})
we have that $0<\pi_1(\Gamma_x|x)\leq f(\Gamma_x|x)\leq
\pi_2(\Gamma_x|x)<1$, which contracts with the nonrandom of $f\in
F$.

(b) By (a)\vspace*{1pt} we only need to show the necessity part. Suppose that $\pi
\in\Pi_s$ and $\eta^\pi\not=\eta^f$ for all
$f\in F$. Then, there exists $D\in{\mathcal{B}}(S)$ such that $0<\hat
\eta^\pi(D)<1$ and $0<\pi(\Gamma_x|x)<1$ for all
$x\in D$ and some $\Gamma_x\in{\mathcal{B}}(A(x))$ (depending on
$x$). Then, by the condition in (b), there exists $x'\in D$ such
that
%
%
\begin{eqnarray}\label{4.21}
0&<&\hat\eta^\pi(\{x'\})<1 \quad\mbox{and}\nonumber\\[-8pt]\\[-8pt]
0&<&\pi(\Gamma_{x'}|x')<1\qquad
\mbox{for some } \Gamma_{x'}\in{\mathcal{B}}(A(x')).\nonumber
\end{eqnarray}
By (\ref{4.21}), we now define two
policies $\pi_1$ and $\pi_2$ as follows:
%
%
\begin{eqnarray}
\label{4.22}
\pi_1(da|x):=\cases{
\pi(da|x), &\quad if $x\not=x'$, \cr
\pi(da\cap\Gamma_{x'}|x')/ \pi(\Gamma_{x'}|x'), &\quad
if $x=x'$;}
\\
\label{4.23}
\pi_2(da|x):=\cases{
\pi(da|x), &\quad if $x\not=x'$, \cr
\pi(da\cap\Gamma_{x'}^c|x')/\pi(\Gamma_{x'}^c|x'), &\quad if $x=x'$.}
\end{eqnarray}

Let $\beta:=\pi(\Gamma_{x'}|x')$, $\delta':=\frac{\beta
\hat\eta^{\pi_2}(\{x'\})}{\beta
\hat\eta^{\pi_2}(\{x'\})+(1-\beta)\hat\eta^{\pi_1}(\{x'\})}$ when
$\hat\eta^{\pi_1}(\{x'\})+\hat\eta^{\pi_1}(\{x'\})>0$, and
$\delta'=\frac{1}{2}$ when
$\hat\eta^{\pi_1}(\{x'\})+\hat\eta^{\pi_1}(\{x'\})=0$. Then, for
each $D\in{\mathcal{B}}(S)$ with $\sup_{x\in D}q^*(x)<\infty$, by
Theorem \ref{Th-4.1} and (\ref{4.22}), (\ref{4.23}) as well as a
straightforward calculation we have
\begin{eqnarray*}
\alpha\hat\eta^{\pi_1}(D)&=&\alpha
\gamma(D)+\int_{S-\{x'\}}q(D|x,\pi)\hat\eta^{\pi_1}(dx)\\
&&{}+\int
_{\Gamma
_{x'}}q(D|x',a)\pi(da|x')\hat\eta^{\pi_1}(\{x'\})/
\beta, \\
\alpha\hat\eta^{\pi_2}(D)&=&\alpha
\gamma(D)+\int_{S-\{x'\}}q(D|x,\pi)\hat\eta^{\pi_2}(dx)\\
&&{}+\int
_{\Gamma
_{x'}^c}q(D|x',a)\pi(da|x')\hat\eta^{\pi_2}(\{x'\})/
(1-\beta).
\end{eqnarray*}
Multiplying by $\delta'$ and $(1-\delta')$ the two equalities,
respectively, and then summarizing, we have
\begin{eqnarray*}
&&\alpha
[\delta'\hat\eta^{\pi_1}(D)+(1-\delta')\hat\eta^{\pi
_2}(D)]\\
&&\qquad=\alpha
\gamma(D)+\int_Sq(D|x,\pi)[\delta'\hat\eta^{\pi_1}(dx)+(1-\delta
')\hat
\eta^{\pi_2}(dx)],
\end{eqnarray*}
which together with Theorem \ref{Th-4.1}(c) implies $\eta^\pi=\delta
'\eta^{\pi_1}+(1-\delta')\eta^{\pi_2}$. Moreover, by
(\ref{4.21}) we see that $0<\eta^{\pi_1}(\{x'\}\times
\Gamma_{x'})=\hat\eta^{\pi_1}(\{x'\})<1$ and
$\eta^{\pi_2}(\{x'\}\times
\Gamma_{x'})=\hat\eta^{\pi_2}(\{x'\})\pi_2(\Gamma_{x'}|x')=0$.
Hence, $\eta^\pi= \delta'\eta^{\pi_1}+(1-\delta')\eta^{\pi_2}$ is
not an extreme point.
\end{pf*}
\begin{pf*}{Proof of Theorem \ref{Th-5.2}} Let $\phi^*$ be a constrained
optimal policy [by Theorem \ref{Th-5.1}(c)], and
${\mathcal{M}}_o^c(e)$ be the set of all extreme points in
${\mathcal{M}}_o^c$ in (\ref{D}). Since ${\mathcal{M}}_o^c$ has been
proved to be convex compact [by Theorem \ref{Th-4.6}(a) and Lemma
\ref{Th-4.4}]. Thus, by Choquet's theorem \cite{P66},
$\eta^{\phi^*}$ is the barycenter of a probability measure
$\bar\mu$ supported on ${\mathcal{M}}_o^c(e)$. Therefore,
%
%
\begin{equation}\label{m}\qquad
\int_{S\times A}c_0(x,a)\eta^{\phi^*}(dx,da)=\int_{{\mathcal
{M}}_o^c(e)}\biggl(\int_{S\times
A}c_0(x,a)\eta(dx,da)\biggr)\bar\mu(d\eta).
\end{equation}
On the other hand, since $\int_{S\times
A}c_0(x,a)\eta^{\phi^*}(dx,da)\leq\int_{S\times
A}c_0(x,a)\eta(dx,da)$ for all $\eta\in{\mathcal{M}}_o^c(e)$, it
follows from (\ref{m}) that there exists $\eta^*\in
{\mathcal{M}}_o^c(e)$ such that
\[
\int_{S\times A}c_0(x,a)\eta^{\phi^*}(dx,da)= \int_{S\times
A}c_0(x,a)\eta^*(dx,da).
\]
Hence, $\pi^*:=\phi^{\eta^*}$ is also constrained optimal.
Moreover, since $\int_{S\times A}c_n(x,a)\eta(dx$, $da)$ (for each
fixed $1\leq n\leq N$) is linear in $\eta\in{\mathcal{M}}_o$ and thus
can be regarded as a ``hyperplane,'' each extreme point of
${\mathcal{M}}_o^c$ is a convex combination of at most $N+1$ extreme
points in ${\mathcal{M}}_0$. That is, there exists $(N+1)$ numbers
$p_k\geq0$ and stationary policies $f_k\in F$ $(k=1, \ldots,
N+1)$ (using Theorem \ref{Th-3.2}) such that $ \eta^*=p_1
\eta^{f_1}+\cdots+p_{N+1} \eta^{f_{N+1}}, p_1+\cdots+p_{N+1}=1$,
which together with Theorem \ref{Th-3.2} and (\ref{4.2}) completes the
proof of this theorem.
\end{pf*}

%

%
\printaddresses

\end{document}